\documentclass[times]{nlaauth}
\usepackage{amsmath}
\usepackage{amssymb}
\usepackage{amsthm}
\usepackage{moreverb}
\usepackage{subfigure}
\usepackage{bm}
\usepackage{todonotes}
\usepackage{graphicx}
\usepackage[colorlinks,bookmarksopen,bookmarksnumbered,citecolor=red,urlcolor=red]{hyperref}

\newcommand\BibTeX{{\rmfamily B\kern-.05em \textsc{i\kern-.025em b}\kern-.08em
T\kern-.1667em\lower.7ex\hbox{E}\kern-.125emX}}

\def\volumeyear{2010}

\DeclareMathOperator*{\argmin}{arg\,\hspace{-0.6mm}min}

\begin{document}

\runningheads{S. Cools, W. Vanroose}{Local Fourier Analysis for Complex Shifted Laplacian}

\title{Local Fourier Analysis of the Complex Shifted Laplacian preconditioner for Helmholtz problems}

\author{Siegfried Cools, Wim Vanroose$^*$}

\address{Department of Mathematics and Computer Science, University of Antwerp, Middelheimlaan 1, 2020 Antwerp, Belgium\vspace{-0.8cm}}

\corraddr{(e-mail) siegfried.cools@ua.ac.be, wim.vanroose@ua.ac.be}


\begin{abstract}
  In this paper we solve the Helmholtz equation with multigrid
  preconditioned Krylov subspace methods. The class of Shifted
  Laplacian preconditioners are known to significantly speed-up Krylov
  convergence. However, these preconditioners have a parameter $\beta \in
    \mathbb{R}$, a measure of the complex shift. Due to
  contradictory requirements for the multigrid and Krylov convergence,
  the choice of this shift parameter can be a bottleneck in applying
  the method. In this paper, we propose a wavenumber-dependent
  minimal complex shift parameter which is predicted by a rigorous
  $k$-grid Local Fourier Analysis (LFA) of the multigrid scheme.  We claim that, given
  any (regionally constant) wavenumber, this minimal complex shift parameter provides the
  reader with a parameter choice that leads to efficient Krylov
  convergence.  Numerical experiments in one and two spatial
  dimensions validate the theoretical results. It appears that the
  proposed complex shift is both the minimal requirement for a
  multigrid V-cycle to converge, as well as being near-optimal in terms of Krylov
  iteration count.
\end{abstract}


\maketitle

\section{Introduction}
The propagation of waves through a material can be mathematically
modelled by the Helmholtz equation
\begin{equation*} 
\left(- \Delta -k(x)^2\right) u(x) = f(x),
\end{equation*}
on a $d$-dimensional subdomain $\Omega \subset \mathbb{R}^d$. Here
$k(x)$ is the space dependent wave number that models change in
material properties as a function of $x\in \Omega$.  For high
wavenumbers $k$, the sparse system of linear equations that results
from the discretization of this PDE is distinctly indefinite, causing most of the
classic iterative methods to perform poorly. Incited by its
efficiency on positive definite problems, the more advanced multigrid
method has repeatedly received attention as possible solver for these
systems. However, as stated in \cite{elman2002multigrid} and
\cite{brandt1986multigrid}, a direct application of the multigrid
method to the discretized systems will inevitably result in
divergence, as both the smoother and the coarse grid correction tends
to introduce growing components in the error.

The aim of this paper is to understand the required modifications 
to the multigrid method such that it effectively solves the indefinite linear
systems resulting from the discretization of the Helmholtz equation.

Over the past few years, many different methods have been proposed to solve
this non-trivial problem, an overview of which can be found in
\cite{ernst2012difficult}. Krylov subspace methods like GMRES
\cite{saad1986gmres} or BiCGStab \cite{van1992bicg} are known for
their ability to definitely converge to the solution. However, Krylov
methods are generally not competitive without a good
preconditioner. Consequently, in recent literature a variety of
specifically Helmholtz-tailored preconditioners have been proposed to
speed up Krylov convergence. In this paper we focus on the class of
so-called Shifted Laplacian preconditioners, which were introduced in
\cite{erlangga2004class} and further analyzed in
\cite{erlangga2006novel} and \cite{erlangga2006comparison}, where they
are shown to be very efficient Krylov method preconditioners,
significantly reducing the number of Krylov iterations on both
academic and realistic problems. Furthermore, unlike the original
Helmholtz system, the Shifted Laplacian preconditioning system can be
solved effectively by a multigrid method.  

Similar to many other Helmholtz solvers like e.g.~\cite{haber2011fast}
however, the Shifted Laplacian method introduces a parameter which is
intrinsic to the method itself, namely the magnitude of the (real
and/or imaginary) shift. Ideally no shift would be introduced, as
preconditioning the original system with its exact inverse causes the
Krylov method to immediately converge to the solution without
iterating. Unfortunately, the use of a (complex) shift is vital to
guarantee multigrid convergence, as stated in
\cite{ernst2012difficult}, \cite{erlangga2004class},
\cite{erlangga2006novel} and \cite{cools2011complex}. Summarizing, one
could say that the Krylov method prefers the shift to be as small as
possible, whereas the multigrid solver only converges for a
sufficiently large shift. This contradiction makes the choice of an
optimal shift a non-trivial yet essential task in fine-tuning the
Shifted Laplacian method.

Note that, whereas in this paper the shift is implemented directly into the discretization matrix, 
it can also be represented by using complex valued grid
distances in the discretization \cite{reps2010indefinite}. Furthermore, a complex
shift can be taken into account in algebraic multilevel method, as elaborated in
\cite{bollhoefer2009algebraic}. The Krylov convergence rate has 
been explored as a function of a complex shift in \cite{osei2010preconditioning}.

The current paper aims at gaining more insight in the value of the
complex shift. up to current date, this shift is determined primeraly on a heuristical
case-by-case basis, as very limited theoretical founding is known.
Furthermore, in recent literature a discussion has risen on the
question whether the choice of this complex shift can be independent
of the wavenumber $k$. We will introduce the notion of a minimal
complex shift, which we define as the smallest possible shift for
multigrid to converge, and show that this shift must indeed depend
directly on the wavenumber. Additionally, it will be verified that the
proposed complex shift is optimal in terms of preconditioned Krylov iterations when
exactly solving the preconditioning problem. When applying a
limited amount of multigrid iterations to the preconditioning problem
as is common practice, only approximately solving the preconditioning
system, it will appear that our definition provides an near-optimal
value for the complex shift parameter w.r.t.~Krylov convergence. 
Consequently, following the given definition, the reader is provided 
a near-optimal and generally safe choice for the complex shift.

The main tool used in this paper to analyze multigrid convergence and
determine a realistic value for the complex shift is the Local Fourier
Analysis, originally introduced by Brandt in 1977
\cite{brandt1977multi}. Our analysis is primarily based on the LFA
setting introduced in \cite{briggs2000multigrid},
\cite{trottenberg2001multigrid} and
\cite{wienands2002three}. Moreover, the analysis in this text can be
seen as an extension of the spectral analysis performed in
\cite{vangijzen2007spectral} and is even more closely related to
recent work like \cite{rodrigo2010accuracy} and
\cite{sheikh2011scalable}, the latter combining Shifted Laplacian with
the new technique of multigrid deflation for improved convergence.

The results in this paper are limited to problems with homogeneous
Dirichlet boundary conditions. Realistic Helmholtz problems often 
have absorbing boundary conditions implemented with an absorbing 
layer such as PML \cite{berenger1994perfectly}. Note that these 
absorbing boundary conditions often lead to a more favorable 
spectrum for Krylov convergence. The analysis performed here
can thus be considered a worst-case convergence scenario for multigrid preconditioned Krylov methods.

This paper is organized as follows. In Section 2, we present a
rigorous LFA analysis of the problem, ultimately allowing us to define
and effectively calculate the minimal complex shift parameter. In
Section 3, numerical results for actually solving the Helmholtz model
problem with a variety of different wavenumbers and grid sizes are
presented, which serve as a validation of the definition given in
Section 2. These experiments verify that the theoretical complex shift
is indeed minimal w.r.t.~multigrid convergence, and we furthermore
observe near-optimality w.r.t.~Krylov convergence. Additionally,
we briefly discuss the similarity between our results and the
observations made in \cite{reps2012analyzing}. We conclude the paper
with a short review of the conclusions in Section 4.

Note that throughout the text, we will use the symbol $\iota$ to denote the 
imaginary unity or $\sqrt{-1}$, to avoid confusion with the index designator $i$.

\section{Local Fourier Analysis}
\subsection{Problem statement and notation}
The general aim of our research is to solve  the $d$-dimensional indefinite Helmholtz equation on an open bounded domain $\Omega \subset \mathbb{R}^d$
\begin{equation} \label{eq:mod1}
- \Delta u + \sigma u = f \quad ~\textrm{on}~ \Omega,
\end{equation}
where $\sigma \in \mathbb{R}^-_0$ is a distinctly negative constant. Observe that we denote the squared wavenumber as $-\sigma = k^2$.\footnote{To avoid unnecessary terminological complications, in this text we loosely refer to $\sigma$ as `the wavenumber'. However, the reader should keep in mind that we hereby intrinsically designate the negatively signed squared wavenumber $-k^2$.} We do not attend to boundary conditions, as Local Fourier Analysis does not take into account the domain boundaries, see Section \ref{sec:C}. Using multigrid on a Complex Shifted Laplacian preconditioner, we consider the problem of iteratively solving the related $d$-dimensional complex shifted Helmholtz equation
\begin{equation} \label{eq:mod2}
- \Delta u + \tilde{\sigma} u = g \quad ~\textrm{on}~ \Omega,
\end{equation}
where $\tilde{\sigma} = \sigma (1+\beta \iota) \in \mathbb{C}$, with complex shift parameter $\beta \in \mathbb{R}^{+}$. Note that the CSL preconditioner was originally introduced in \cite{erlangga2004class} with a somewhat more general $\tilde{\sigma} = \sigma (\alpha+\beta \iota)$, but following the observations in  \cite{erlangga2006novel} we will choose to permanently set $\alpha \equiv 1$. This is in a way a natural choice for $\alpha$, as preserving the real shift keeps the preconditioner very close to the original problem. Equation (\ref{eq:mod2}) is typically discretized using a finite difference scheme, yielding a linear system $A \bold{u} = \bold{g}$. In our multigrid analysis, we assume a standard second-order finite difference discretization on an equidistant mesh with $N$ gridpoints, where typically $N = 2^p$ for some $p \in \mathbb{N}_0$ and with mesh width $h = 1/N$. This renders for $d=1$ a discretization matrix $A^{1D}$ with stencil representation
\begin{equation} \label{eq:1Dsten}
A^{1D} = \frac{1}{h^2} \left( \begin{array}{ccc}   -1 & 2+\tilde{\sigma} h^2 & -1   \end{array} \right),
\end{equation}
or, for $d=2$, a discretization matrix $A^{2D}$ given by
\begin{equation} \label{eq:2Dsten}
A^{2D} = \frac{1}{h^2} \left( \begin{array}{ccc}  & -1 & \\ -1 & 4+\tilde{\sigma} h^2 & -1 \\ & -1 &   \end{array} \right) ,
\end{equation}
which can easily be generalized for higher dimensions. Throughout the analysis, we consider the propagation of the initial fine grid error $\bold{e}_1^{(0)}$ through a $k$-grid analysis error propagation
matrix $M^k_1$ as follows
\begin{equation} \label{eq:eprop}
\bold{e}_1^{(m+1)} = M^k_1~\bold{e}_1^{(m)}, \qquad m \geq 0,
\end{equation}
where $M^k_l$ designates the iteration matrix on grid level $l$ corresponding to the multigrid cycle that employs $k$ grid levels. Note that we designate the finest grid by the index $l=1$, the one-level coarser grid by $l=2$, and so on, implying a total of $2^{p-l+1}$ grid points on the $l$-th grid. For the well-known two-grid cycle, $M_1^2$ can be expressed as (see \cite{briggs2000multigrid})
\begin{equation} \label{eq:propM2}
M_1^2 = S_1^{\nu_2} (I_1 - I_2^1 A_2^{-1} I_1^2 A_1) S_1^{\nu_1},
\end{equation}
which is generalized by the
following definition of the error propagation matrix $M_l^k$ for a
$k$-grid analysis (see
\cite{trottenberg2001multigrid},\cite{wienands2002three},\cite{rodrigo2010accuracy})
\begin{align} \label{eq:propM}
M_l^k &= S_l^{\nu_2} (I_l - I_{l+1}^l (I_{l+1} - M_{l+1}^k)A_{l+1}^{-1} I_l^{l+1} A_l) S_l^{\nu_1},\qquad l = 1,\ldots,k-1, \notag \\
M_k^k &= 0,
\end{align}
where $S_l$ is the smoothing operator, $A_l$ is the discretization matrix representation and $I_l$ is the identity matrix on the $l$-th coarsest grid $G_l$ with mesh size $h_l = 2^{l-1}h ~ (l=1,\ldots,k)$. For the one-dimensional problem, the $l$-th coarsest grid $G_l$ is given by 
\begin{equation} \label{eq:gridl}
G_l = \{x \in \Omega ~|~ x = x_j = j h_l, ~ j \in \mathbb{Z}\}, 
\end{equation}
whereas in 2D it is defined as
\begin{equation} \label{eq:gridl2}
G_l = \{\bold{x}  \in \Omega ~|~ \bold{x} = (x_{j_1}, x_{j_2}) = (j_1 h_l, j_2 h_l), ~ j_1,j_2 \in \mathbb{Z}\}.
\end{equation}
Furthermore, $I_l^{l+1}: G_{l} \to G_{l+1}$ and $I_{l+1}^l: G_{l+1} \to G_{l}$ are the restriction and prolongation operators, respectively. In this paper, we will use full weighting restriction and linear interpolation as the standard intergrid operators. 

\subsection{Basic principles of 1D Local Fourier Analysis} \label{sec:C}

We briefly sketch the key ideas behind Local Fourier Analysis. LFA is based on the assumption that both relaxation and two-grid correction are local processes, in which each unknown is updated using only the information in neighbouring points. Furthermore, boundary conditions are neglected by extending all multigrid components to an infinite grid. It is presumed that the error $\bold{e}_{l}^{(m)}$ on the $l$-th coarsest grid can be written as a formal linear combination of the Fourier modes $\varphi_l(\theta,x) = e^{\iota \theta x/h_l}$ with $x \in G_l$ and Fourier frequencies $\theta \in \mathbb{R}$. These frequencies may be restricted to the interval $\Theta = (-\pi,\pi] \subset \mathbb{R}$ as a consequence of the fact that for $x \in G_l$
\begin{equation} \label{eq:ffreq1}
\varphi_l(\theta+2\pi,x) = \varphi_l(\theta,x).
\end{equation}
The set of $l$-th grid Fourier modes is typically denoted
\begin{equation} \label{eq:space}
\mathcal{E}_l = \text{span}\{\varphi_l(\theta,x) = e^{\iota \theta x/h_l} \, | \, x\in G_l, \, \theta \in \Theta\}.
\end{equation}
The Fourier modes are known \cite{stüben1982multigrid} to be formal eigenfunctions of the operator $A_l$. More precisely, the general relation $A_l \varphi_l(\theta,x) = \tilde{A}_l(\theta) \varphi_l(\theta,x)$ holds, where the formal eigenvalue $\tilde{A}_l(\theta)$ is called the Fourier symbol of the operator $A_l$. Given a so-called low Fourier frequency $\theta^0 \in (-\pi/2,\pi/2]$, its complementary frequency $\theta^1$ is defined as
\begin{equation} \label{eq:theta1}
\theta^1 = \theta^0 - \text{sign}(\theta^0) \pi.
\end{equation}
Using this notation, the following important property can easily be derived (for $1 \leq l < k$)
\begin{equation} \label{eq:ffreq2}
\varphi_l(\theta^0,x) = \varphi_{l+1}(2\theta^0,x) = \varphi_l(\theta^1,x) , \qquad x \in G_{l+1}.
\end{equation}
The Fourier modes $\varphi_l(\theta^0,x)$ and $\varphi_l(\theta^1,x)$ are called $(l+1)$-th level harmonic modes. These Fourier modes coincide on the $(l+1)$-th coarsest grid, where they are represented by a single coarse grid mode with double frequency $\varphi_{l+1}(2\theta^0,x)$. In this way, each low-frequency mode $\varphi_l(\theta^0,x)$ is naturally associated with a high-frequency mode $\varphi_l(\theta^1,x)$ on the $l$-th grid. It is convenient to denote the two-dimensional subspace of $\mathcal{E}_l$ spanned by these $(l+1)$-th level harmonics as
\begin{equation} \label{eq:subsp}
\mathcal{E}_l^{\theta^0} = \text{span}\{\varphi_l(\theta^0,\cdot),\varphi_l(\theta^1,\cdot)\}, \qquad 1 \leq l < k,
\end{equation}
with $\theta^1$ as defined by (\ref{eq:theta1}). Note that every coarse-level subspace $\mathcal{E}_{l+1}^{\theta^0}$ can be decomposed into spaces spanned by finer-level harmonics, as (\ref{eq:ffreq2}) implies 
\begin{equation} \label{eq:subsp2}
\mathcal{E}_{l+1}^{\theta^0} = \mathcal{E}_{l}^{\theta^0/2} \cup \mathcal{E}_{l}^{\theta^1/2}, \qquad 1 \leq l < k-1.
\end{equation}
The significance of these spaces $\mathcal{E}_l^{\theta}$ is that they are invariant under both smoothing operators and correction schemes under general assumptions. Due to the previous observations, throughout the Fourier symbol calculation one can assume without loss of generality that each $l$-th grid error $\bold{e}_l^{(m)}$ can be decomposed into components $e_{l,j}^{(m)} = \bold{e}^{(m)}_l(x_j)$ that consist of a single Fourier mode $\varphi_l(\theta,x_j)$ and thus can be represented as
\begin{equation} \label{eq:ecomf}
e_{l,j}^{(m)} = A^{(m)} \, \varphi_l(\theta,x_j) = A^{(m)} \, e^{\iota j \theta}, \qquad \theta \in \Theta, \quad j \in \mathbb{Z}, \quad m \geq 0,
\end{equation}
where the amplitude $A^{(m)}$ changes as a function of the iteration $m$. Moreover, using this expression together with the error relation in (\ref{eq:eprop}) and the stencil representation of the operators forming the error propagation matrix $M_1^k$, it can be derived that the error amplitudes in two subsequent iterations are related by an amplification factor
function $\mathcal{G}_l(\theta,\sigma,\beta)$
\begin{equation} \label{eq:Gampl}
A^{(m+1)} = \mathcal{G}_l(\theta,\sigma,\beta) A^{(m)}, \qquad m \geq 0,
\end{equation}
which describes the evolution of the amplitude $A^{(m)}$ through consecutive iterations. It is our aim to elaborate an analytical expression for the amplification factor on the finest grid $\mathcal{G}_1(\theta,\sigma,\beta)$, which can be considered a continuation of the spectrum of $M_1^k$ (see \cite{briggs2000multigrid}, \cite{cools2011complex}). For notational convenience, we tend to denote $\mathcal{G}_1(\theta,\sigma,\beta)$ as $\mathcal{G}(\theta,\sigma,\beta)$, dropping the fine grid index. This amplification factor can be computed by calculating the Fourier symbols of each component of the $G_1$-grid error propagation matrix $M_1^k$. After the symbols are structured in so-called eigenmatrices for each component, combining harmonic frequencies, and the appropriate products of these matrices are taken, a $2^{k-1}\times2^{k-1}$ eigenmatrix representation $\tilde{M}_1^k$ of the 1D error propagation matrix $M_1^k$ is obtained. The amplification factor function $\mathcal{G}_1(\theta,\sigma,\beta)$ can then readily be computed as the spectral radius of $\tilde{M}_1^k$ (see \cite{trottenberg2001multigrid},\cite{stüben1982multigrid})
\begin{equation} \label{eq:Gexpr}
\mathcal{G}(\theta,\sigma,\beta) = \rho(\tilde{M}_1^k) = \max|\lambda(\tilde{M}_1^k)|.
\end{equation}
A well-known condition for a general iterative method with given iteration matrix, say $M$, to have convergence, is that $\rho(M) < 1$. This convergence condition can be generalised to the following demand on the $G_1$-grid amplification factor
\begin{equation} \label{eq:Gcond}
\max_{\theta \in \Theta} \, \mathcal{G}(\theta,\sigma,\beta) \leq 1,
\end{equation}
clearly representing the analogue condition for a multigrid cycle. We denote the amplication factor as a function of the Fourier frequency $\theta$, the wavenumber $\sigma$ and the complex shift parameter $\beta$.

\textbf{Definition of the minimal shift.} With the help of the amplification factor we define the \emph{minimal complex shift  parameter} $\beta_{\min}(\sigma)$ as a function of the wavenumber
$\sigma$ as
\begin{equation} \label{eq:defbmin}
\beta_{\min} := \argmin_{\beta \geq 0}\left\{\max_{\theta \in \Theta} \, \mathcal{G}(\theta,\sigma,\beta) \leq 1\right\}.
\end{equation}
This definition can be interpreted through (\ref{eq:Gcond}) as the smallest possible complex shift required for the multigrid method to converge, i.e.~it is the smallest value of $\beta$ for which every
single eigenmode of the error is reduced through consecutive multigrid iterations. Additionally, the numerical experiments presented in Section 3 will show that the complex shift parameter $\beta_{\min}$ as
defined here is near-optimal for any multigrid-preconditioned Krylov method in terms of iteration count. This means that when the preconditioner is inverted exactly (up to discretization error) using a sufficiently large amount of multigrid steps, a minimal number of Krylov steps is required when choosing the value of the complex shift equal to $\beta_{\min}$. 

\subsection{The 1D Fourier symbols} \label{sec:A}

In this section we will effectively calculate the Fourier symbols of the different component operators of a 1D $k$-grid scheme.\\

\textbf{Discretization operator.} The evolution of the error under the discretization operator $A_l$ can be calculated from its stencil representation (\ref{eq:1Dsten}), yielding
\begin{equation*} \label{eq:discre}
e_{l,j}^{(m+1)} = \frac{1}{h_l^2}\left(-e_{l,j-1}^{(m)}+(2+\tilde\sigma h_l^2) \, e_{l,j}^{(m)}-e_{l,j+1}^{(m)}\right), \qquad j \in \mathbb{Z}.
\end{equation*}
Using expression (\ref{eq:ecomf}), the amplitude relation is found to be 
\begin{equation*}
A^{(m+1)} = \left(-\frac{2}{h_l^2}\cos\theta+\frac{2}{h_l^2}+\tilde\sigma\right) A^{(m)},
\end{equation*}
and hence the discretization operator Fourier symbol $\tilde{A}_l(\theta)$ is
\begin{equation} \label{eq:discrG}
\tilde{A}_l(\theta) = \left(-\frac{2}{h_l^2}\cos\theta+\frac{2}{h_l^2}+\tilde\sigma\right).
\end{equation}

\textbf{Restriction operator.} Using an analogous stencil argument, the error propagation under application of the full weighting restriction operator $I_l^{l+1}$ can be derived to be
\begin{equation*} \label{eq:restre}
e_{l+1,j}^{(m+1)} = \frac{1}{4} \, e_{l,2j-1}^{(m)} + \frac{1}{2} \, e_{l,2j}^{(m)} + \frac{1}{4} \, e_{l,2j+1}^{(m)}, \qquad j \in \mathbb{Z}.
\end{equation*}
Substituting the error components by (\ref{eq:ecomf}), one obtains the Fourier symbol $\tilde{I}_l^{l+1}(\theta)$
\begin{equation} \label{eq:restrG}
\tilde{I}_l^{l+1}(\theta) = \frac{1}{2}(\cos\theta+1).
\end{equation}

\textbf{Interpolation operator.} The linear interpolation operator $I_{l+1}^l$ propagates the coarse grid error as
\begin{equation*} \label{eq:intere}
e_{l,2j}^{(m+1)} = e_{l+1,j}^{(m)} \quad \textrm{or} \quad e_{l,2j+1}^{(m+1)} = \frac{1}{2} e_{l+1,j}^{(m)} + \frac{1}{2} e_{l+1,j+1}^{(m)} \, , \qquad j \in \mathbb{Z},
\end{equation*} 
where the first equation holds if $x \in G_l \cap G_{l+1}$ and the second if $x \in G_l \backslash G_{l+1}$. Once again substituting the error components by (\ref{eq:ecomf}) and combining the above cases, the following expression for the Fourier symbol $\tilde{I}_{l+1}^l(\theta)$ is obtained
\begin{equation} \label{eq:interG}
\tilde{I}_{l+1}^l(\theta) = \frac{1}{2}(\cos\theta+1).
\end{equation}
Note that full-weighting restriction and linear interpolation are dual operators, yielding the exact same Fourier symbol $\tilde{I}_l^{l+1}(\theta)=\tilde{I}_{l+1}^l(\theta)$.\\

\textbf{Smoothing operator.} As a smoother we will use standard weighted Jacobi relaxation throughout this text. It is easily derived (see e.g.~\cite{briggs2000multigrid}) that the $\omega$-Jacobi relaxation matrix $S_l = I_l-\omega D_l^{-1}A_l$ relates the error in a gridpoint $x_j \in G_l$ through subsequent iterations by
\begin{equation*} \label{eq:jace}
e_{l,j}^{(m+1)} = (1-\omega) e_{l,j}^{(m)} + \frac{\omega}{2+\tilde\sigma h_l^2} (e_{l,j-1}^{(m)} + e_{l,j+1}^{(m)}), \qquad j \in \mathbb{Z}.
\end{equation*}
Presuming the error is of the form (\ref{eq:ecomf}), one obtains the following amplitude relation and $\omega$-Jacobi smoother Fourier symbol $\tilde{S}_l(\theta)$
\begin{equation} \label{eq:jacG}
\tilde{S}_l(\theta) = \left(1-\omega+\frac{2\omega}{2+\tilde\sigma h_l^2}\cos\theta\right).
\end{equation}
For completeness and generality, we note that the following analysis can readily be conducted using a different smoother scheme like Gauss-Seidel relaxation, for which an analogous calculation shows the Fourier symbol $\tilde{S}_l(\theta)$ is given by
\begin{equation} \label{eq:gssG}
\tilde{S}_l(\theta) = \left(\frac{e^{i\theta}}{2+\tilde\sigma h_l^2 - e^{-i\theta}}\right).
\end{equation}

Once calculated, the 1D Fourier symbols of each component of $M_1^k$ can now be structured in eigenmatrices, representing the action of a $l$-th grid component on the subspace of $(l+1)$-grid harmonics, where $l = 1,\ldots,k-1$. This matrix-building will be performed explicitely in the next sections for specific values of $k$ to compute the $G_1$-grid amplification function $\mathcal{G}(\theta,\sigma,\beta)$ and, from this function, the minimal complex shift parameter $\beta_{\min}$.

\subsection{A 1D Local Fourier Analysis of the weighted Jacobi smoother}

\begin{figure}[t] \centering
\includegraphics[width=6cm]{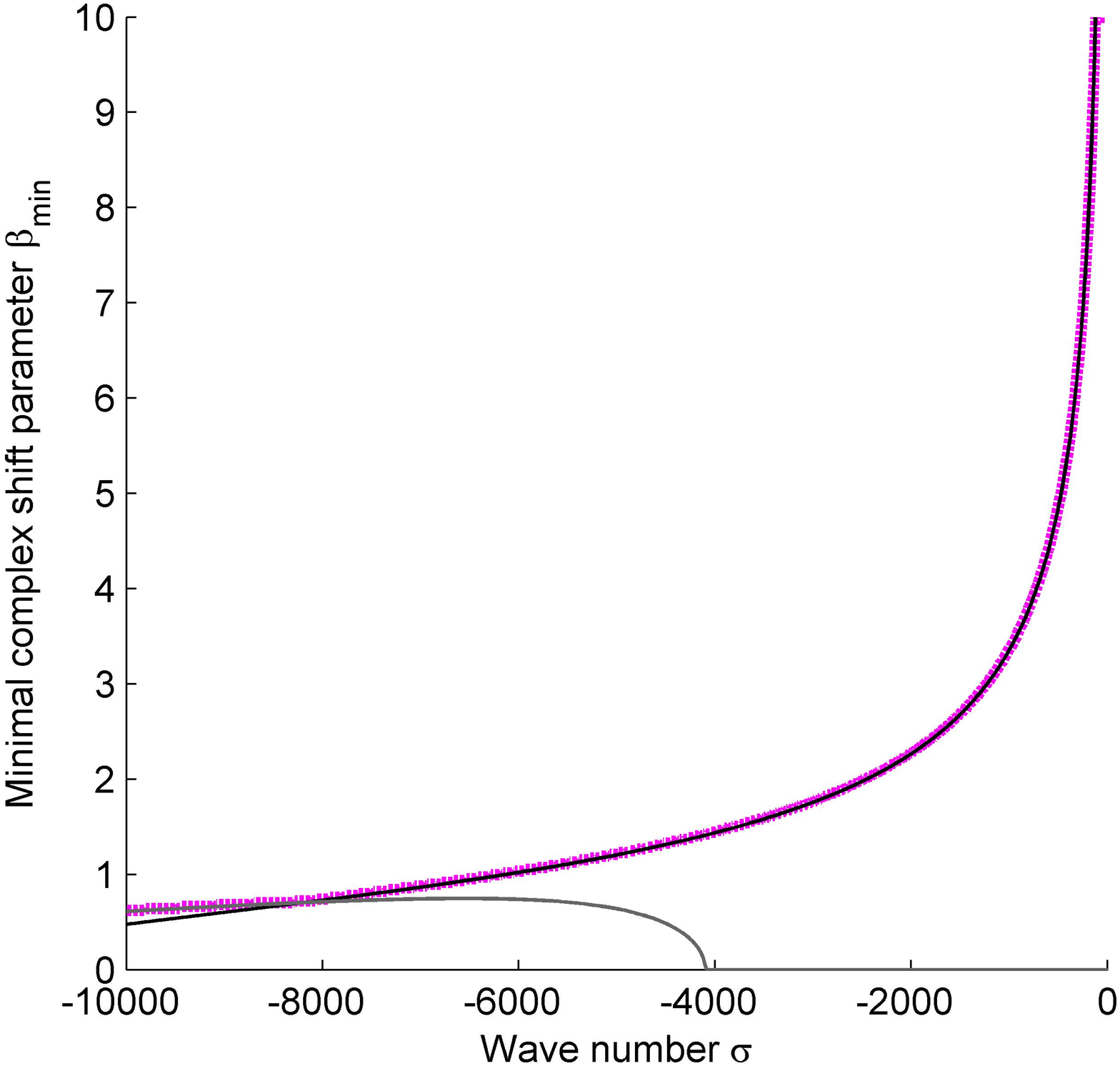}
\captionsetup{width=14.1cm}
\caption[The minimal complexe shift parameter for 1D smoother]{\textsl{The minimal complex shift parameter $\beta_{\min}$ from the 1D smoother analysis with $N=64$ and $\omega = 2/3$ as a function of the wavenumber $\sigma$. Solid black line: analytical lower limit (\ref{eq:lowlim}) corresponding to Fourier frequency $\theta = 0$. Solid grey line: analytical lower limit (\ref{eq:lowlim2}) corresponding to the frequency $\theta = \pi$.}}
\label{fig:1Dbetasmoother}
\end{figure}

To get some insight in the problem, we will initially perform a very basic Local Fourier Analysis of the $\omega$-Jacobi smoother. It is well-known that the smoother, forming an essential part of any multigrid scheme, is unstable for the indefinite Helmholtz problem. This instability is caused by the smoothest eigenmodes which have negative eigenvalues and consequently diverge under the action of the smoother, as described to some extend in \cite{elman2002multigrid} and \cite{brandt1986multigrid}. To make up for smoother instability, one could determine a complex shift based purely on the divergence of the eigenmodes under application of the fine-grid smoother operator $S_1$. Hence, recalling that every two-dimensional subspace of coarse grid harmonics $\mathcal{E}_1^{\theta^0} \subset \mathcal{E}$ with $\theta^0 \in (-\pi/2,\pi/2]$ is left invariant under the actions of the smoother operator $S_1$, the fine grid error propagation eigenmatrix can be written as
\begin{equation} \label{eq:jemtrx}
\left[ \begin{matrix} \tilde{S}_1(\theta^0) & 0 \\ 0 & \tilde{S}_1(\theta^1) \end{matrix} \right],
\end{equation}
with $\tilde{S}_1(\theta)$ being the fine grid smoother symbol derived in (\ref{eq:jacG}). Since the matrix is diagonal, the calculation of the spectral radius in (\ref{eq:Gexpr}) is easy. This can then
be substituted in definition (\ref{eq:defbmin}) to obtain the minimal complex shift $\beta_{\min}$. Choosing this value for $\beta$, one ensures that all eigenmodes converge under the action of $S_1$. The result is shown in Figure \ref{fig:1Dbetasmoother}, where the minimal complex shift parameter $\beta_{\min}$ is plotted as a function of the wave number $\sigma$ for $N = 64$ gridpoints.

One observes from Figure \ref{fig:1Dbetasmoother} that in order to  compensate for smoother instability, values of $\beta_{\min}$ are generally very large with $\lim_{\sigma \to 0} \beta_{\min} =
+\infty$, suggesting extremely large complex shifts are needed to eliminate smoother instability. Due to the simple form of the eigenmatrix (\ref{eq:jemtrx}), the lower limit $\beta_{\min}$ can easily
be calculated analytically. The frequency $\theta$ maximizing $|\tilde{S}_1(\theta)|$ over $(-\pi,\pi]$ can be found by setting the first derivative of $|\tilde{S}_1(\theta)|$ equal to zero, rendering
\begin{equation*}
-8\omega^2\cos\theta\sin\theta - 4\omega(1-\omega)(2+\sigma h^2)\sin\theta = 0,
\end{equation*}
which reveals $\theta = 0, \, \theta = \pi$ and $\theta = \pm\arccos\{-(1-\omega)(2+\sigma h^2)/2 \omega\}$ as local extrema. A second derivative check confirms $\theta = 0$ and $\theta = \pi$ to maximize the smoother symbol $|\tilde{S}_1(\theta)|$. The value of $\beta_{\min}$ can be derived by substituting both maxima in inequality (\ref{eq:Gcond}). Presuming $\theta = 0$ the expression for $\mathcal{G}$ reduces to 
\begin{equation*}
\mathcal{G}(0,\sigma,\beta) = \left|1-\omega+\frac{2\omega}{2+\tilde\sigma h^2}\right| \leq 1,
\end{equation*}
which can be elaborated using the definition $\tilde\sigma = \sigma(1+\beta \iota)$ to yield a lower limit on the complex shift
\begin{equation} \label{eq:lowlim}
\sqrt{\frac{4}{(\omega-2)\sigma h^2} -1} \leq \beta.
\end{equation}
A similar lower limit can be derived for the maximum in $\theta = \pi$, which eventually comes down to
\begin{equation} \label{eq:lowlim2}
\sqrt{-\frac{\omega(4+\sigma h^2)^2 -2(2+\sigma h^2)(4+\sigma h^2)}{(\omega-2)(\sigma h^2)^2}} \leq \beta.
\end{equation}
The minimal complex shift $\beta_{\min}$ is now defined as the minimal value of $\beta$ satisfying both inequalities. Note that from the left-hand of inequality (\ref{eq:lowlim}), it readily follows that $\lim_{\sigma \to 0} \beta_{\min} = +\infty$, as we observed from Figure \ref{fig:1Dbetasmoother}.

Although this short analysis of the smoother is useful to obtain some intuition, the value of the minimal complex shift $\beta_{\min}$ is clearly overestimated since we only consider the smoother operator, while in a multigrid setting the smoothest error components are removed by the coarse grid correction. Indeed, to guarantee multigrid convergence a much smaller value of $\beta$ may be chosen compared to
the shift suggested by Figure \ref{fig:1Dbetasmoother}. In the next sections the second component of the multigrid method, the correction scheme, will be taken into account to attain a realistic
curve for $\beta_{\min}$ as a function of the wave number.

\subsection{A 1D two-grid Local Fourier Analysis} \label{sec:B}
A more realistic curve for the minimal complex shift of a multigrid cycle is obtained by considering a basic 2-grid analysis with only one presmoothing step, for which the fine grid error propagation matrix is given by 
\begin{equation} \label{eq:2grpro}
M_1^2 = (I_1 - I_2^1 A_2^{-1} I_1^2 A_1) S_1.
\end{equation}
As suggested earlier, every two-dimensional subspace of coarse grid harmonics $\mathcal{E}_1^{\theta^0} \subset \mathcal{E}$ with $\theta^0 \in (-\pi/2,\pi/2]$ is left invariant under the actions of both the smoother operator $S_1$ and the pure 2-grid correction operator $I_2 - I_2^1 A_2^{-1} I_1^2 A_1$. The action of the total 2-grid error propagation matrix $M_1^2$ on $\mathcal{E}_1^{\theta^0}$ is given by its $2\times2$ eigenmatrix
\begin{equation} \label{eq:2greig}
\tilde{M}_1^2 = \left[
I_1 -
\left[ \begin{matrix} \tilde{I}_2^1(\theta^0) \\ \tilde{I}_2^1(\theta^1) \end{matrix} \right] 
\tilde{A}_2(2\theta^0)^{-1}
\left[ \begin{matrix} \tilde{I}_1^2(\theta^0) \\ \tilde{I}_1^2(\theta^1) \end{matrix} \right]^T
\left[ \begin{matrix} \tilde{A}_1(\theta^0) \\ \tilde{A}_1(\theta^1) \end{matrix} \right]^D
\right] 
\left[ \begin{matrix} \tilde{S}_1(\theta^0) \\ \tilde{S}_1(\theta^1) \end{matrix} \right]^D.
\end{equation}
We use the superscript-$^D$ notation to designate the operation that transforms a vector into a diagonal matrix by placing the entries of the superscribed vector along the main diagonal. The spectral radius of this expression, and thus the function $\mathcal{G}(\theta,\sigma,\beta)$, can easily be calculated analytically or numerically.

\begin{figure}[t] \centering
\subfigure[]{\includegraphics[width=3.5cm]{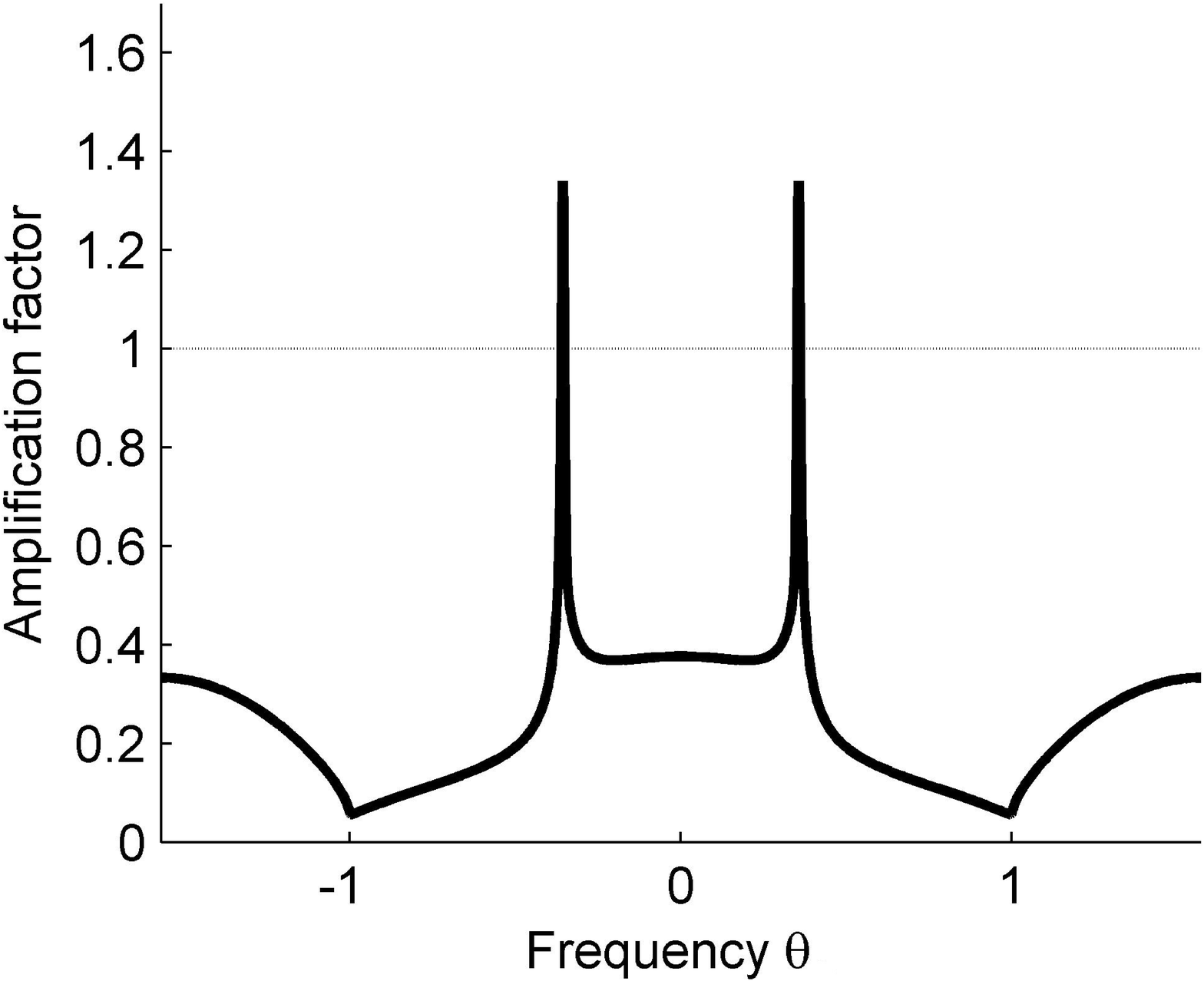}}
\subfigure[]{\includegraphics[width=3.5cm]{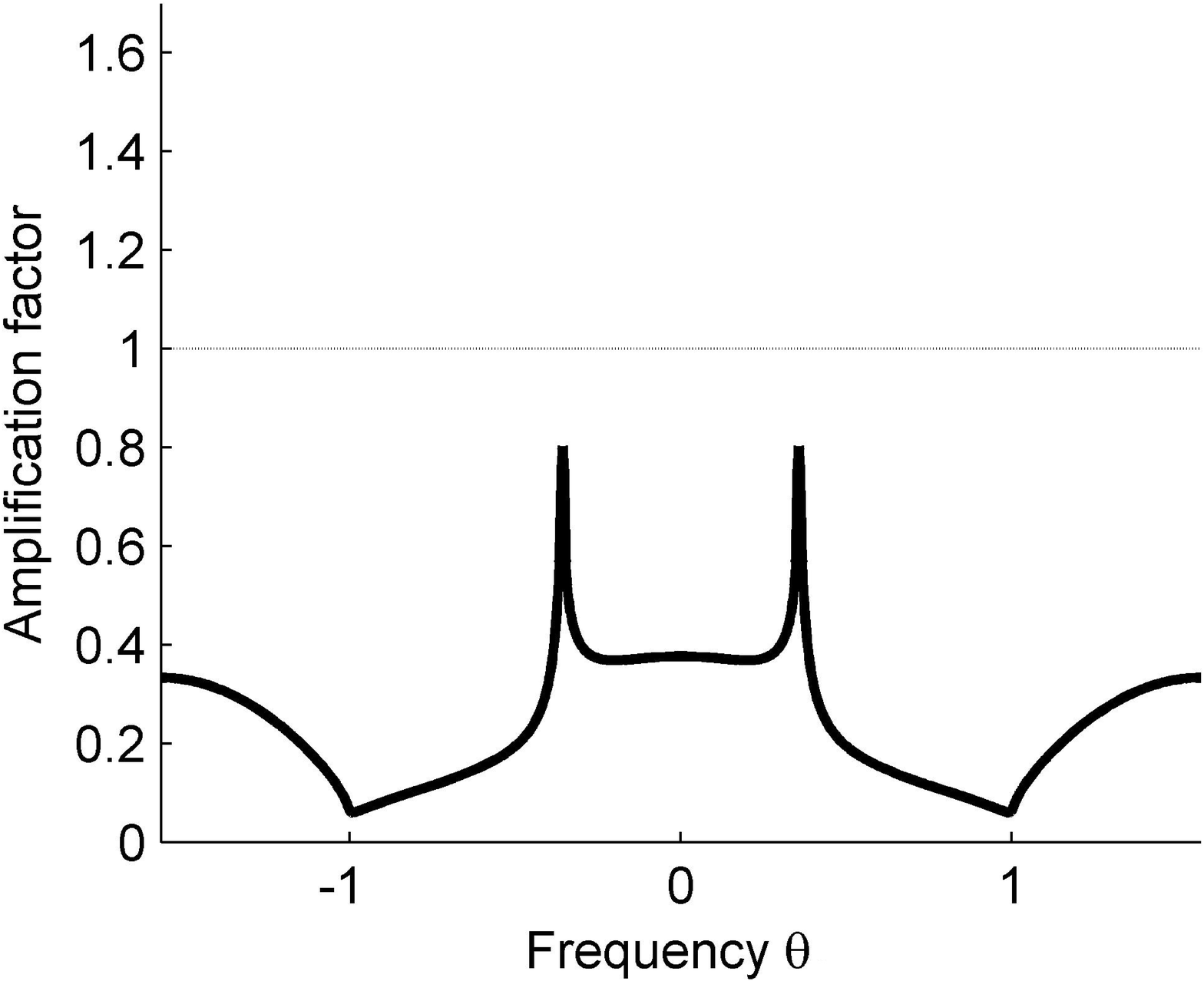}}
\subfigure[]{\includegraphics[width=3.5cm]{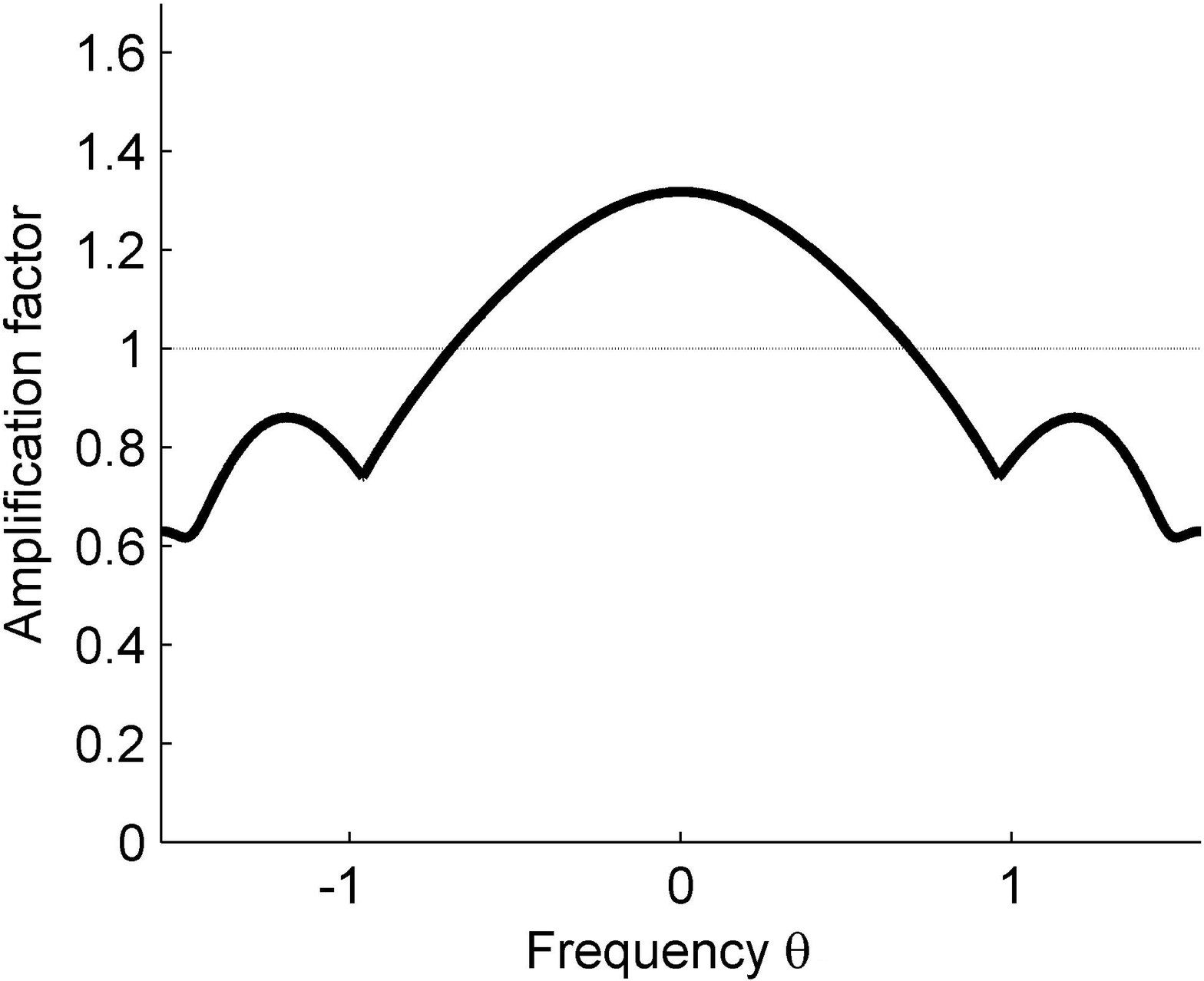}}
\subfigure[]{\includegraphics[width=3.5cm]{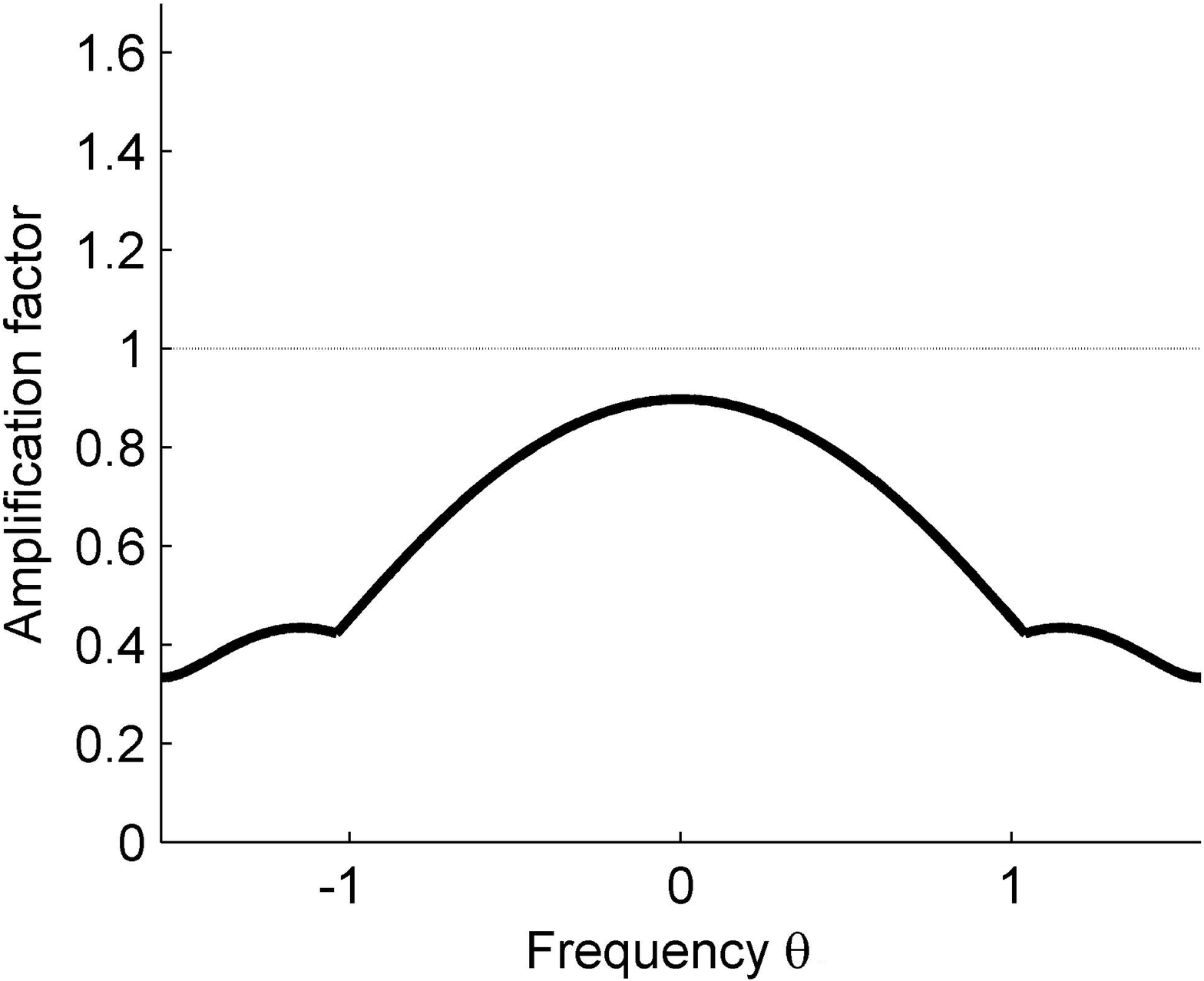}}
\captionsetup{width=14.1cm}
\caption[1D 2-grid \mathcal{G}(\theta,\sigma,\beta)]{\textsl{Amplification factor function $G(\theta,\sigma,\beta)$ for the 1D two-grid analysis with $N=64$ as a function of $\theta \in (-\pi/2,\pi/2]$ for fixed wavenumber and complex shift. Two leftmost figures: wavenumber $\sigma = -500$ and complex shift parameter $\beta = 0.02$ (a) and $\beta = 0.04$ (b). Two rightmost figures: $\sigma = -5000$ and $\beta = 0.2$ (c) and $\beta = 0.8$ (d).}}
\label{fig:plotG}
\end{figure}

In Figure \ref{fig:plotG} we show the amplification factor for two choices of $\sigma$. In panel (a) and (b) of the figure, where $\sigma=-500$, we notice the appearance of a resonance that is caused by the coarse grid correction. Only this one resonant mode tends to diverge, whereas for the majority of the frequencies $\theta$ the amplification factor is significantly smaller than 1. The appearance of such a resonance was already discussed in \cite{ernst2012difficult} and originates from the inversion of the coarse grid discretization symbol $\tilde{A}_2(2\theta^0)$ in (\ref{eq:2greig}). Using expression (\ref{eq:discrG}), the frequency $\theta$ corresponding to the semi-asymptote can be approximated by 
\begin{equation} \label{eq:approx1} 
\theta \approx \pm \arcsin\sqrt{-\sigma h_2^2/4} = \pm \arcsin\sqrt{-\sigma h^2},
\end{equation}
forcing the real part of the denominator of $\mathcal{G}$, i.e.~the real part of the symbol $\tilde{A}_2(2\theta^0)$, to zero and thus maximizing the function $\mathcal{G}(\theta,\sigma,\beta)$ for fixed
values of $\sigma$ and $\beta$. Note that approximation (\ref{eq:approx1}) only makes sense when $|\sigma| < 1/h^2$. This discussion suggests that for small values of $|\sigma|$ the divergence of the coarse grid correction scheme is to be the determining factor in the choice of the complex shift $\beta_{\min}$.

Alternatively, for high values of $|\sigma|$ the maximum of
$G(\theta,\sigma,\beta)$ appears to be spread broadly around $\theta =
0$, as one perceives from Figure \ref{fig:plotG} (c) and (d), which
implies a large range of smooth eigenmodes tends to diverge. The
maximum here is not due to the coarse grid correction, but originates
primarily from the divergence of the smoothing operator. To
substantiate this statement, we consider the smoothing operator
Fourier symbol $\tilde{S}_1(\theta)$ calculated in
(\ref{eq:jacG}). Presuming $-4/h^2 < \sigma < -1/h^2$, approximation
(\ref{eq:approx1}) no longer holds (no semi-asymptotes are generated),
and it can easily be established using derivative arguments that
without a complex shift (i.e.~with $\tilde\sigma = \sigma$),
expression (\ref{eq:jacG}) reaches a maximum in either $\theta = 0$
(if $\sigma > -2/h^2$) or its complementary frequency $\theta = \pi$
(if $\sigma < -2/h^2$). Note that the maximization of
$\tilde{S}_1(\theta)$ was already discussed in the previous
section. Additionally, one can verify that the absolute value of this
maximum is always larger than 1, implying divergence of the smoothest
mode when applying the smoothing operator $S_1$. It is clear that for
large values of $|\sigma|$, the main factor determining the value of
$\beta_{\min}$ is the divergence caused by the smoother operator,
rather than the coarse grid correction.

\begin{figure}[t] \centering
\includegraphics[width=6cm]{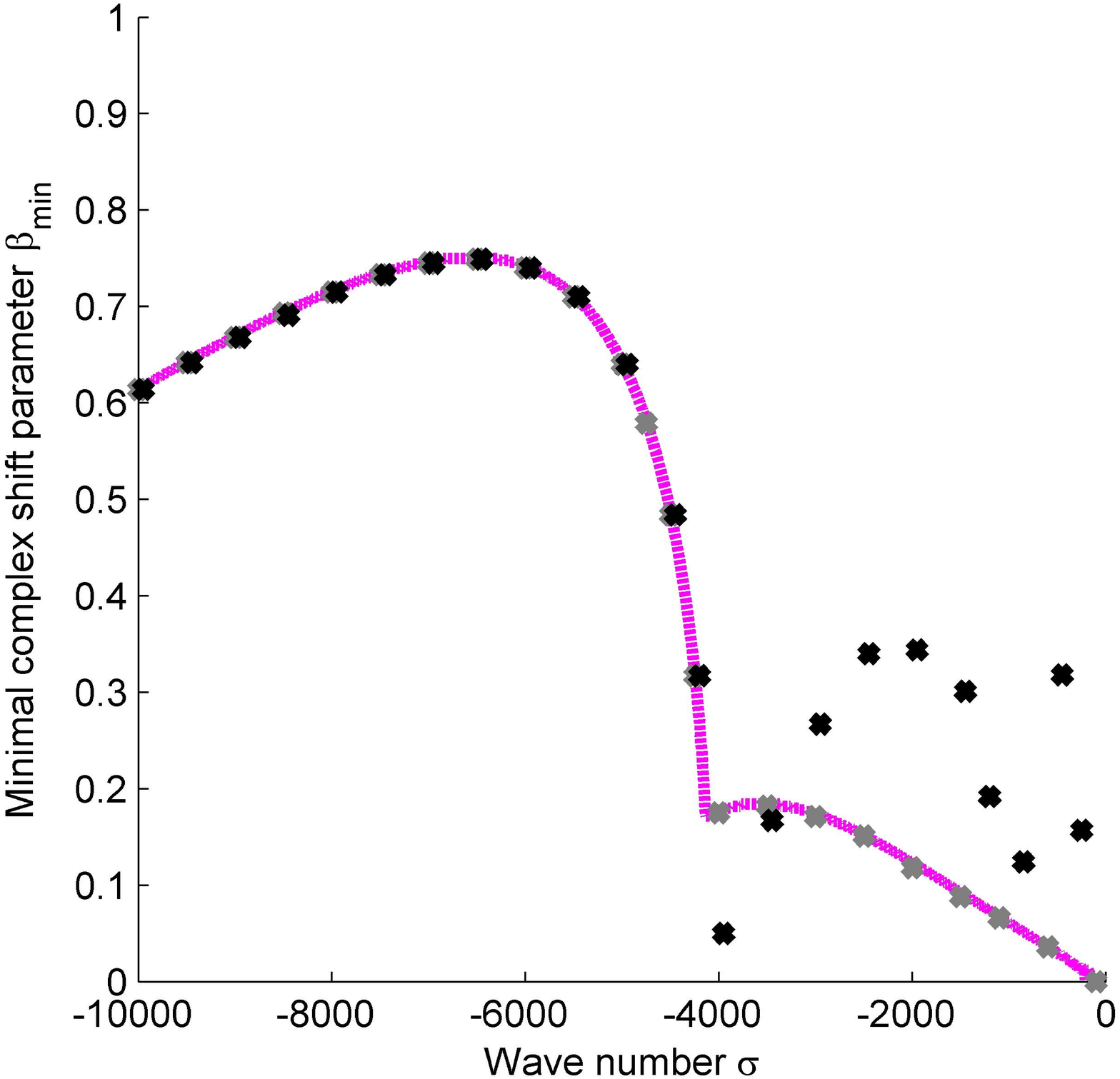}
\captionsetup{width=14.1cm}
\caption[The minimal complexe shift parameter for 1D 2-grid]{\textsl{The minimal complex shift parameter $\beta_{\min}$ from the 1D two-grid analysis with $N=64$ and $\omega = 2/3$ as a function of the wavenumber $\sigma$. Black and grey points represent experimentaly measured values of $\beta_{\min}$ (see text).}}
\label{fig:1Dbeta2grid}
\end{figure}

From analytical expression (\ref{eq:2greig}) for the amplification factor $\mathcal{G}(\theta,\sigma,\beta)$, one can now numerically calculate the value of $\beta_{\min}$ given a fixed wavenumber $\sigma$. The maximization over $\theta \in \Theta$ is done using an equidistant discretization of the frequency domain with sufficiently small frequency-step. Subsequently, the minimal complex shift argument $\beta_{\min}$ is computed using 10 steps of the elementary bisection method, providing an accuracy of at least 3 decimals to $\beta_{\min}$. Figure \ref{fig:1Dbeta2grid} represents the minimal complex shift parameter $\beta_{\min}$ as a function of the wavenumber $\sigma$ for $N = 64$ gridpoints. Notice the significant increase of the curve around $\sigma = -1/h^2$, corresponding to the difference in the underlying component generating $\beta_{\min}$ as discussed above: for wavenumbers $\sigma < -1/h^2$ the divergence of the smoothing operator is the determining factor in selecting $\beta_{\min}$\footnote{Note the striking resemblance between the high wavenumber regime of Figure \ref{fig:1Dbeta2grid} and the analytical lower limit (\ref{eq:lowlim2}) corresponding to $\theta = \pi$ displayed on Figure \ref{fig:1Dbetasmoother} (solid grey line).}, whereas for $\sigma > -1/h^2$ the application of the coarse grid correction operator is decisive.

Furthermore, we added to the figure some results for $\beta_{\min}$ based on the experimentally measured convergence rate. The gray points in overplot represent
these experimentally measured values of $\beta_{\min}$,
determined by numerical calculation of the asymptotic convergence
factors for the two-grid scheme on a random initial fine-grid error
$\bold{e}_1^{(0)}$, whilst subjecting those factors to criterion
(\ref{eq:Gcond}). The experimental values are perfectly matched by the
theoretical curve. The black points represent some similarly computed
experimental values for a full V-cycle. 

It is clear from Figure \ref{fig:1Dbeta2grid} that the aforementioned analysis yields accurate results for the two-grid
scheme. However, it is also clear that the addition of multiple coarser grids completely alters the choice of $\beta_{\min}$ for 
small values of $|\sigma|$. To obtain a sufficiently accurate simulation of the full
V-cycle, a higher order $k$-grid scheme should be applied. In
particular, we would like to guarantee a near-to-exact theoretical
prediction of $\beta_{\min}$ around $\sigma = -(0.625/h)^2$ for the
V-cycle, respecting the $kh\leq 0.625$ criterion from
\cite{bayliss1985accuracy} for a minimum of 10 gridpoints per
wavelength on the standard $[0,1]$-domain. From this perspective, a
4-grid LFA analysis will appear to be satisfactory (see next section). We
emphasize however that this paper primarily focusses on the iterative
solution, rather than the accuracy of the discretization. Indeed, the high wavenumber regime $\sigma < -4000$ of the $\beta_{\min}$ curve shown by Figure \ref{fig:1Dbeta2grid} is primarily useful from a theoretical point of view. However it must be noted that when descending to coarser levels in the multigrid hierarchy, wavenumbers look relatively larger compared to the finer levels.

\subsection{A 1D 3-grid and 4-grid Local Fourier Analysis}

\begin{figure}[t] \centering
\includegraphics[width=6cm]{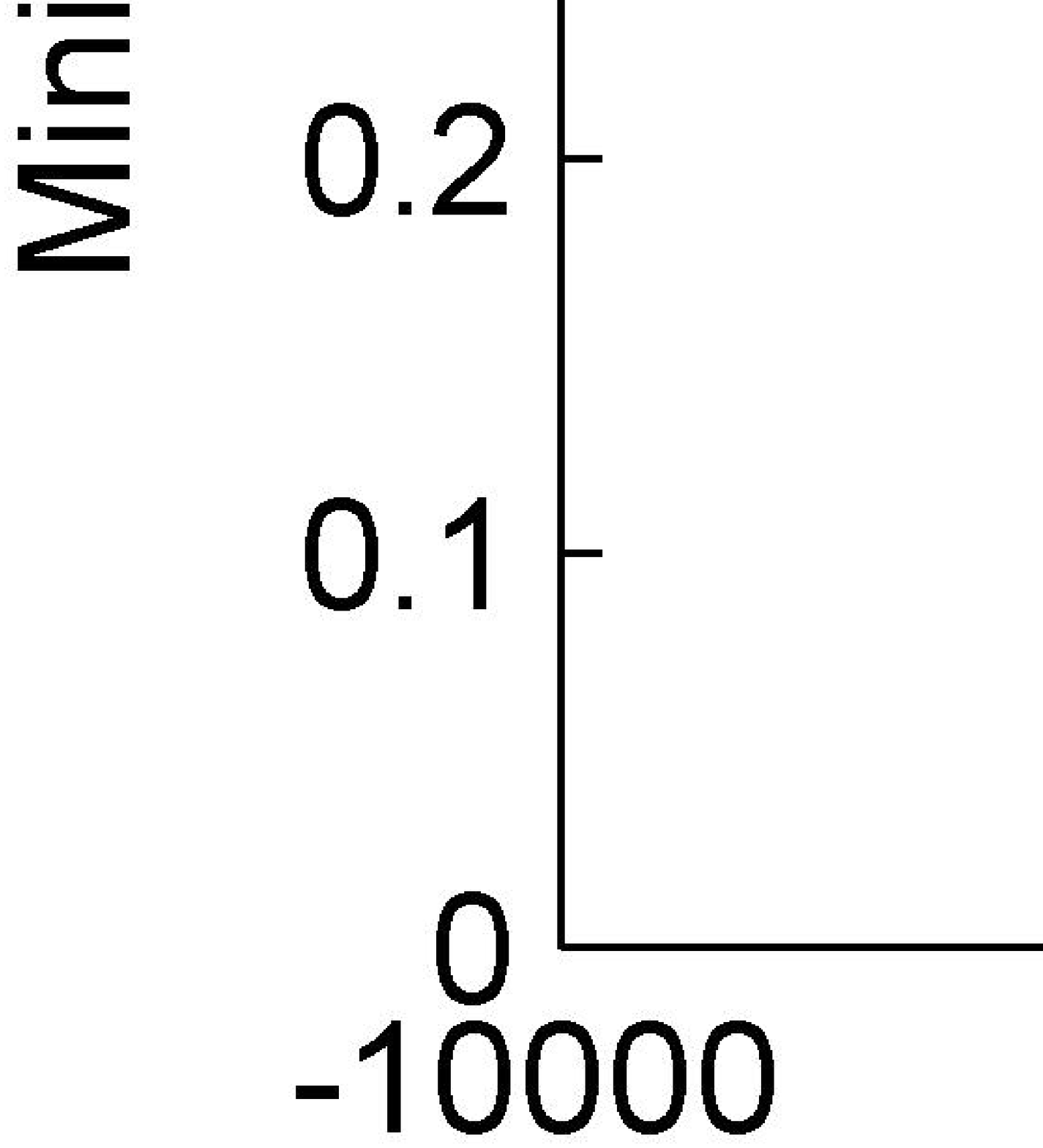} \quad
\includegraphics[width=6cm]{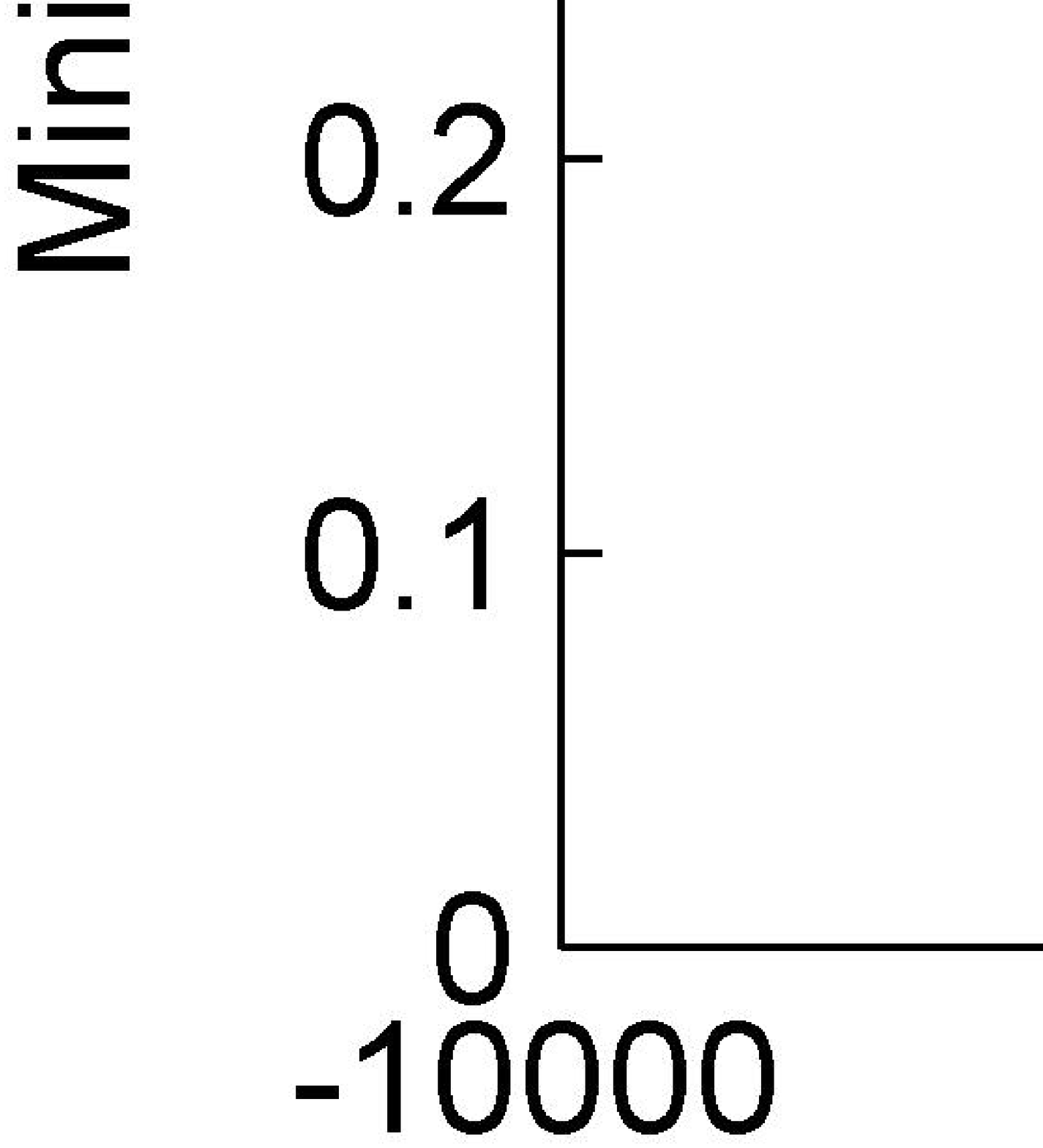}
\captionsetup{width=14.1cm}
\caption[The minimal complex shift parameter for 1D 3-grid and 4-grid]{\textsl{The minimal complex shift parameter $\beta_{\min}$ from the 1D 3-grid (left) and 4-grid (right) analysis with $N=64$ and $\omega = 2/3$ as a function of the wavenumber $\sigma$.}}
\label{fig:1Dbeta34grid}
\end{figure}

In this section we extend the two-grid LFA analysis to a more general $k$-grid analysis, motivated by the inaccurate results of a two-grid analysis for low values of $|\sigma|$. For notational purposes, we restrict ourselves to a rigourous explanation of the $k=3$ case, however the 4-grid analysis is completely analogous. Results are shown in Figure \ref{fig:1Dbeta34grid} for both $k=3$ and $k=4$.

The 3-grid fine grid error propagation matrix is given by (\ref{eq:propM}) to be
\begin{equation} \label{eq:3grpro}
M_1^3 = (I_1 - I_{2}^1 (I_2 - (I_2 - I_{3}^2 A_{3}^{-1} I_2^3 A_2) S_2)A_{2}^{-1} I_1^{2} A_1) S_1.
\end{equation}
For any $\theta^0 \in (-\pi/2,\pi/2]$, the 3-grid operator leaves the 4-dimensional subspace $\mathcal{E}_2^{\theta^0} \subset \mathcal{E}$ of $G_3$ harmonics invariant, implying its eigenmatrix $\tilde{M}_1^3$ is a $4\times4$ square matrix
\begin{align} \label{eq:3greig}
\tilde{M}_1^3 &=
\left[
I_1 -
\left[ \begin{matrix} \tilde{I}_2^1(\theta^{0,0}) ~~ 0 \\ \tilde{I}_2^1(\theta^{0,1}) ~~ 0 \\ 0 ~~ \tilde{I}_2^1(\theta^{1,0}) \\ 0 ~~ \tilde{I}_2^1(\theta^{1,1}) \end{matrix} \right] 
\left[
I_2 - \tilde{M}_2^3
\right]
{\left[ \begin{matrix} \tilde{A}_2(\theta^{0}) \\ \tilde{A}_2(\theta^{1}) \end{matrix} \right]^D}^{-1}
\left[ \begin{matrix} \tilde{I}_2^1(\theta^{0,0}) ~~ 0 \\ \tilde{I}_2^1(\theta^{0,1}) ~~ 0 \\ 0 ~~ \tilde{I}_2^1(\theta^{1,0}) \\ 0 ~~ \tilde{I}_2^1(\theta^{1,1}) \end{matrix} \right]^T
\left[ \begin{matrix} \tilde{A}_1(\theta^{0,0}) \\ \tilde{A}_1(\theta^{0,1}) \\ \tilde{A}_1(\theta^{1,0}) \\ \tilde{A}_1(\theta^{1,1}) \end{matrix} \right]^D
\right]
\left[ \begin{matrix} \tilde{S}_1(\theta^{0,0}) \\ \tilde{S}_1(\theta^{0,1}) \\ \tilde{S}_1(\theta^{1,0}) \\ \tilde{S}_1(\theta^{1,1}) \end{matrix} \right]^D \hspace{-0.2cm}, \notag
\\
\tilde{M}_2^3 &=
\left[
I_2 -
\left[ \begin{matrix} \tilde{I}_3^2(\theta^0) \\ \tilde{I}_3^2(\theta^1) \end{matrix} \right] 
\tilde{A}_3(2\theta^0)^{-1}
\left[ \begin{matrix} \tilde{I}_2^3(\theta^0) \\ \tilde{I}_2^3(\theta^1) \end{matrix} \right]^T
\left[ \begin{matrix} \tilde{A}_2(\theta^0) \\ \tilde{A}_2(\theta^1) \end{matrix} \right]^D
\right] 
\left[ \begin{matrix} \tilde{S}_2(\theta^0) \\ \tilde{S}_2(\theta^1) \end{matrix} \right]^D \hspace{-0.2cm},
\end{align}
where we introduce as a short notation $\theta^{0,0} = \theta^0/2$ and $\theta^{1,0} = \theta^1/2$, and additionally $\theta^{0,1} = \theta^0/2+\pi$ and $\theta^{1,1} = \theta^1/2+\pi$. The amplification factor function $\mathcal{G}(\theta,\sigma,\beta)$ is per definition the spectral radius of the eigenmatrix $M_1^3$. 

We remark that the resonances in the function $\mathcal{G}$ discussed in the previous section indeed reappear in the three-level approximation. However, approximating the resonant frequency analogue to (\ref{eq:approx1}) becomes a significantly non-trivial matter for $k>2$ due to the inversion of multiple Fourier symbols.

Given a wavenumber $\sigma$, the value of the corresponding minimal complex shift $\beta_{\min}$ can then be calculated numerically as described in the previous section. The theoretical results from the 3-grid and 4-grid analysis are shown on Figure \ref{fig:1Dbeta34grid} for $N=64$. The black dots again represent experimentally measured values of $\beta_{\min}$ for a full V-cycle. Note that the theoretical curve from the 4-grid analysis precisely matches the experimental values around $\sigma = -(0.625/h)^2 = -1500$.

From Figure \ref{fig:1Dbeta34grid} one observes that adding an additional level in the Local Fourier Analysis reveals a new `kink' coupled by a local minimum in the $\beta_{\min}$ curve. Note how every additional local minimum revealed by the $k$-grid operator $M_1^k$ is situated at the exact value of $\sigma$ for which the coarsest grid operator $A_k$ turns completely negative definite. Indeed, the 2-grid operator displays a local minimum around $\sigma = -4096$, the 3-grid analysis reveals an additional local minimum around $\sigma = -1024$, and the 4-grid analysis adds a minimum at $\sigma = -256$. Taking into account more levels in the $k$-grid analysis would indeed result in further refinement of the $\beta_{\min}$ curve in the low-wavenumber regime. However, we have restricted our analysis to the 4-grid scheme, as is custom in the literature.

\subsection{Basic principles of 2D Local Fourier Analysis}

The Local Fourier Analysis performed in the previous sections can easily be extended to 2D problems. In this case the error $\bold{e}_l^{(m)}$ can be written as a formal linear combination of the 2D Fourier modes $\varphi_l(\bm{\theta},\bold{x}) = e^{\iota (\theta_1 x_{j_1} + \theta_2 x_{j_2})/h_l}$ with $\bold{x} = (x_{j_1}, x_{j_2})\in G_l$, which is now defined by (\ref{eq:gridl2}), and a couple of frequencies ${\bm\theta} = (\theta_1,\theta_2) \in \Theta^2 = (-\pi,\pi]^2 \subset \mathbb{R}^2$. In analogy to (\ref{eq:space}), the subspace of $l$-th grid Fourier modes is denoted
\begin{equation} \label{eq:space2}
\mathcal{E}_l = \text{span}\{\varphi_l(\bold{x},\bm{\theta}) = e^{\iota (\theta_1 x_{j_1} + \theta_2 x_{j_2})/h_l} \, | \, \bold{x} =(x_{j_1}, x_{j_2}) \in G_l, \, \bm{\theta}= (\theta_1,\theta_2) \in \Theta^2\}.
\end{equation}
Considering a low frequency $\bm{\theta}^{00} = (\theta^{00}_1,\theta^{00}_2) \in (-\pi/2,\pi/2]^2$ and defining its 2D complementary frequencies as
\begin{align}
\bm{\theta}^{10} &= \bm{\theta}^{00} - (\text{sign}(\theta^{00}_1) \pi, 0), \notag \\
\bm{\theta}^{01} &= \bm{\theta}^{00} - (0, \text{sign}(\theta^{00}_2) \pi), \\
\bm{\theta}^{11} &= \bm{\theta}^{00} - (\text{sign}(\theta^{00}_1) \pi, \text{sign}(\theta^{00}_2 \pi)), \notag
\end{align}
one can derive the harmonics property for $1 \leq l < k$
\begin{equation} \label{eq:harm2}
\varphi_l(\bm{\theta}^{00},\bold{x}) = \varphi_l(\bm{\theta}^{10},\bold{x}) = \varphi_{l+1}(2\bm{\theta}^{00},\bold{x}) = \varphi_l(\bm{\theta}^{01},\bold{x}) = \varphi_l(\bm{\theta}^{11},\bold{x}), \qquad \bold{x} \in G_{l+1}.
\end{equation}
The low frequency Fourier mode $\varphi_l(\bm{\theta}^{00},\bold{x})$ and three high frequency modes $\varphi_l(\bm{\theta}^{10},\bold{x})$, $\varphi_l(\bm{\theta}^{01},\bold{x})$ and $\varphi_l(\bm{\theta}^{11},\bold{x})$ coincide on the $(l+1)$-th coarsest grid and thus are called $(l+1)$-th level harmonics. Analogous to (\ref{eq:subsp}), we denote the four-dimensional subspace of $\mathcal{E}_l$ spanned by these $(l+1)$-th level harmonics as 
\begin{equation} \label{eq:subsp3}
\mathcal{E}_l^{\bm{\theta}^{00}} = \text{span}\{\varphi_l(\bm{\theta}^{00},\cdot),\varphi_l(\bm{\theta}^{10},\cdot),\varphi_l(\bm{\theta}^{01},\cdot),\varphi_l(\bm{\theta}^{11},\cdot)\}, \qquad 1 \leq l < k.
\end{equation}
Again, these spaces are invariant under general smoothing operators and correction schemes. Given these notations, the 2D Local Fourier Analysis itself is completely similar to the 1D case, with every $l$-th grid error component $e^{(m)}_{l,(j_1,j_2)}$ being represented as a single Fourier mode
\begin{equation} \label{eq:ecomf2}
e_{l,(j_1,j_2)}^{(m)} = A^{(m)} \, \varphi_l(\bm{\theta},(x_{j_1},x_{j_2})), \qquad \bm{\theta} \in \Theta^2, \quad j_1,j_2 \in \mathbb{Z}, \quad m \geq 0,
\end{equation}
from which a relation between the amplitudes in subsequent iterations can be derived
\begin{equation} \label{eq:Gampl2}
A^{(m+1)} = \mathcal{G}_l(\bm{\theta},\sigma,\beta) A^{(m)}, \qquad m \geq 0.
\end{equation}
The finest grid amplification factor $\mathcal{G}_1(\bm{\theta},\sigma,\beta)$ of the 2D operator $M_1^k$ can be calculated according to (\ref{eq:Gexpr}), where $\tilde{M}_1^k$ now represents the $4^{k-1} \times 4^{k-1}$ eigenmatrix representation of $M_1^k$. The Fourier symbols defining this eigenmatrix are derived in the next section. 

\textbf{Definition of the minimal shift.} Once calculated, the amplification factor function allows us to compute the minimal complex shift parameter, defined in 2D as
\begin{equation} \label{eq:defbmin2}
\beta_{\min} := \argmin_{\beta \geq 0}\left\{\max_{\bm{\theta} \in \Theta^2} \, \mathcal{G}(\bm{\theta},\sigma,\beta) \leq 1\right\},
\end{equation}
which is the analogue of definition (\ref{eq:defbmin}) for the two-dimensional case.

\subsection{The 2D Fourier symbols}

Below we briefly sum the Fourier symbols of the different component operators of the 2D $k$-grid operator $M_1^k$. The elaboration of these symbols is completely similar to the calculations in Section \ref{sec:A}, and is omitted here in favour of readability.\\

\textbf{Discretization operator.} Using stencil representation (\ref{eq:2Dsten}) and expression (\ref{eq:ecomf2}), the discretization operator Fourier symbol $\tilde{A}_l(\bm{\theta})$ can readily be found to be
\begin{equation} \label{eq:discrG2}
\tilde{A}_l(\bm{\theta}) = -\frac{2}{h_l^2}\cos\theta_1-\frac{2}{h_l^2}\cos\theta_2+\frac{4}{h_l^2}+\tilde\sigma.
\end{equation}

\textbf{Restriction operator.} Using an analogous stencil argument, the Fourier symbol of the 2D full weighting restriction operator $I_l^{l+1}$ can be derived as
\begin{equation} \label{eq:restrG2}
\tilde{I}_l^{l+1}(\bm{\theta}) = \frac{1}{4}(\cos\theta_1 \cos\theta_2 + \cos\theta_1 + \cos\theta_2 + 1).
\end{equation}

\textbf{Interpolation operator.} As is the case in the 1D analysis, the Fourier symbol of the 2D linear interpolation operator $I_{l+1}^l$ can be derived to yield exactly the same expression as the full weighting restriction, i.e.
\begin{equation} \label{eq:interG2}
\tilde{I}_{l+1}^l(\bm{\theta}) = \frac{1}{4}(\cos\theta_1 \cos\theta_2 + \cos\theta_1 + \cos\theta_2 + 1).
\end{equation}

\textbf{Smoothing operator.} Presuming the error is of the form (\ref{eq:ecomf2}), one obtains the following $\omega$-Jacobi smoother Fourier symbol for the 2D case
\begin{equation} \label{eq:jacG2}
\tilde{S}_l(\bm{\theta}) = 1-\omega+\frac{2\omega}{4+\tilde\sigma h_l^2}(\cos\theta_1+\cos\theta_2).
\end{equation}

\begin{figure}[t] \centering
\includegraphics[width=6cm]{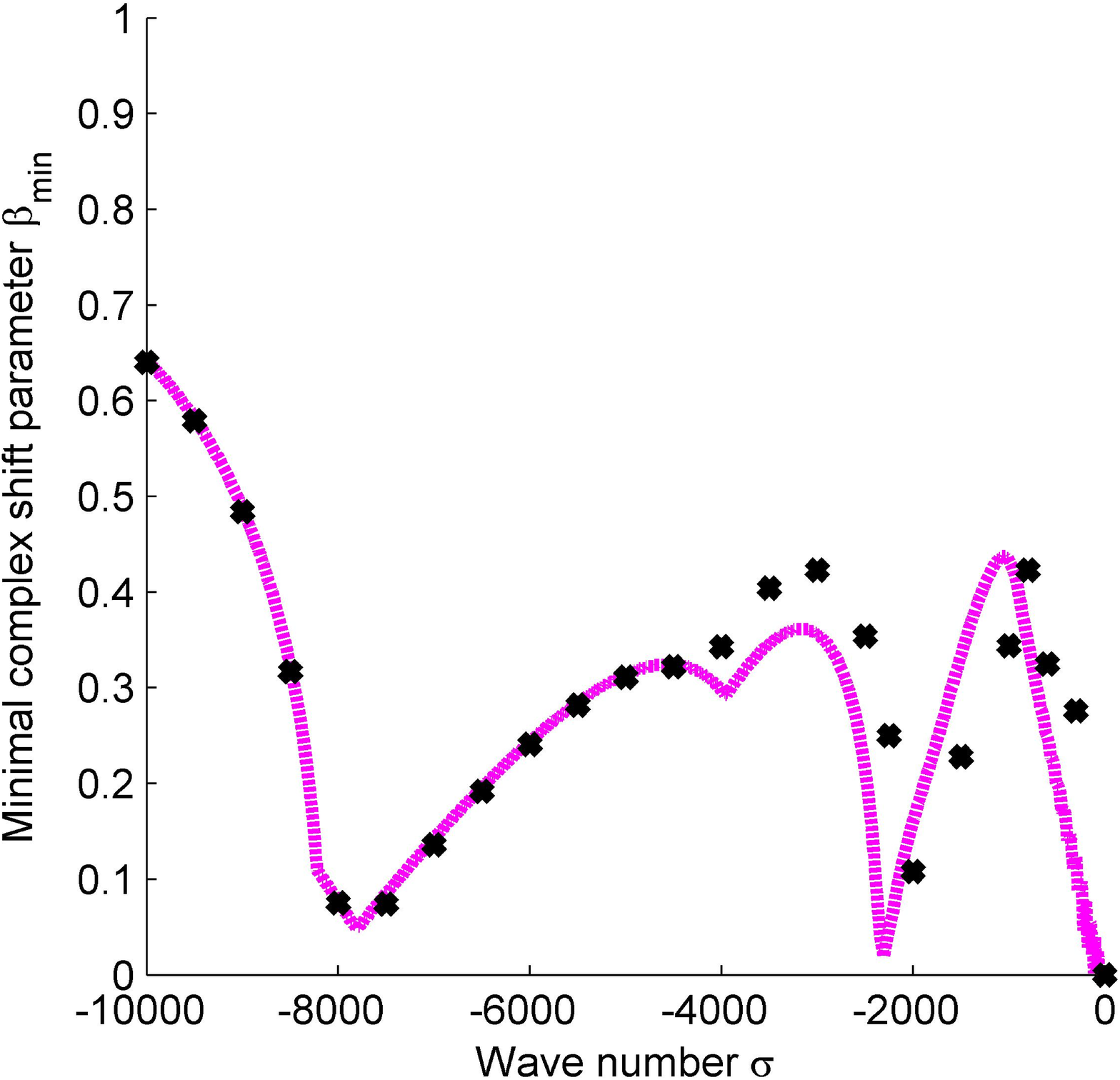} \quad
\includegraphics[width=6cm]{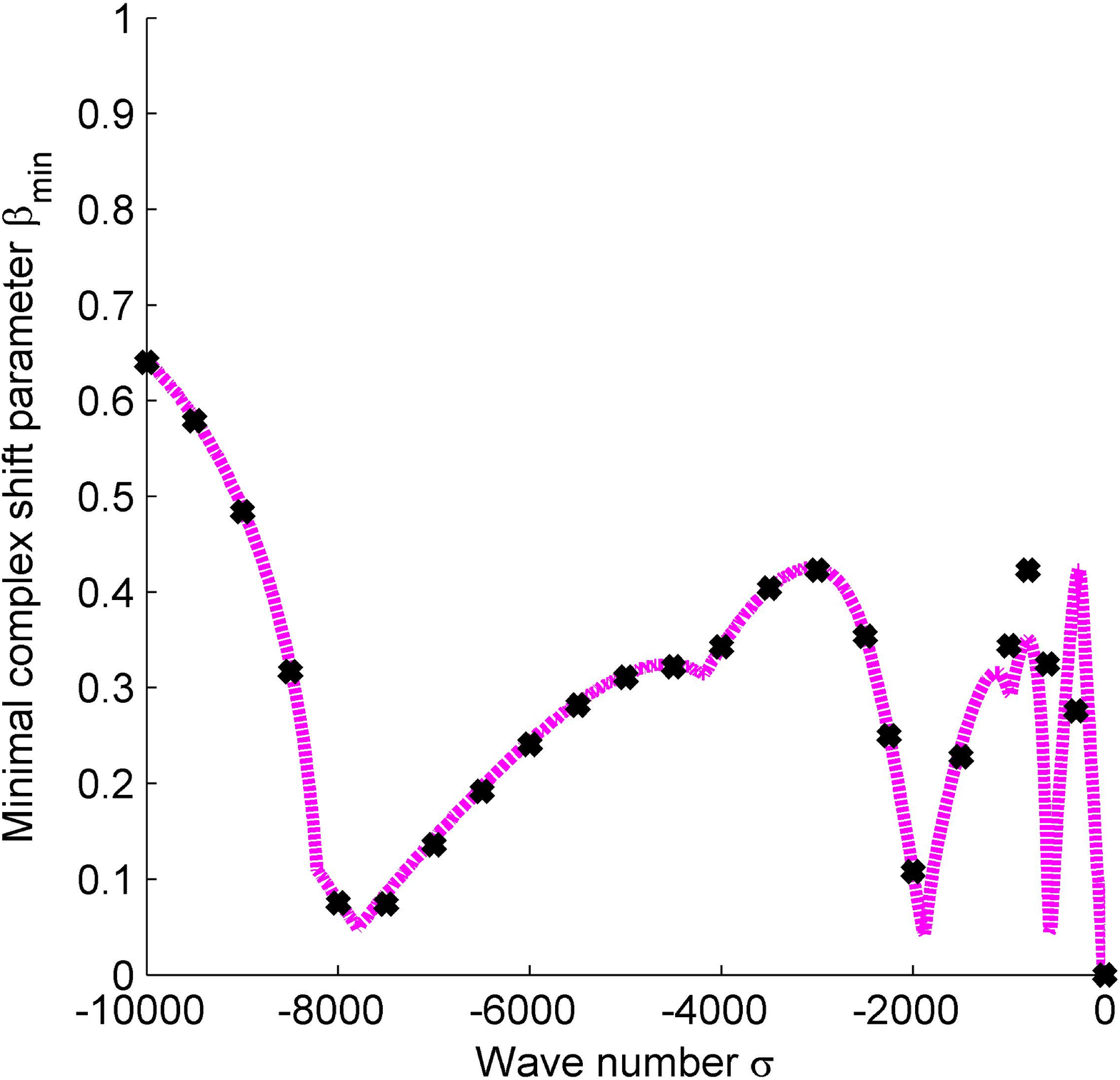}
\captionsetup{width=14.1cm}
\caption[The minimal complexe shift parameter for 2D 3-grid and 4-grid]{\textsl{The minimal complex shift parameter $\beta_{\min}$ from the 2D 3-grid (left) and 4-grid (right) analysis with $N=64$ and $\omega = 2/3$ as a function of the wavenumber $\sigma$.}}
\label{fig:2Dbeta34grid}
\end{figure}

\subsection{A 2D 3-grid and 4-grid Local Fourier Analysis}

Once calculated, the 2D Fourier symbols of each component of $M_1^k$ can again be structured in eigenmatrices, representing the action of the $l$-th grid components on the subspace of $(l+1)$-th level
harmonics, where $l = 1,\ldots,k-1$. The product of these matrices yields the eigenmatrix $\tilde{M}_1^k$, from which the amplification function $\mathcal{G}(\bm{\theta},\sigma,\beta)$ can be calculated. For any given $\sigma$, the value of the minimal complex shift parameter $\beta_{\min}$ can then be computed numerically as described in Section \ref{sec:B}. We have effectively performed the 2D computation for $k = 2,3$ and $4$. The results from the 3-grid and 4-grid analysis for $N=64$ are shown by Figure \ref{fig:2Dbeta34grid}. Clearly, the observations made in the 1D case are again visible here. 

We notice a significant increase of the $\beta_{\min}$ value for wavenumbers $\sigma = -2/h^2$ and smaller. For these wavenumbers, the maximum of the 2D amplification factor is determined solely by the divergence of the smoother operator. For values of $\sigma > -2/h^2$, the resonances caused by the correction scheme that appear in the 1D problem near a single frequency given by equation \eqref{eq:approx1}, now appears for a range of couples $(\theta_1,\theta_2)$ for which the real part of the (combination of) coarse grid symbol(s) is approximately zero. As in the 1D case, a 4-grid analysis accurately simulates $\beta_{\min}$ values extracted from V-cycle experiments up to wavenumbers as small as $\sigma = -1500$.

\subsection{Extensions and general remarks}

\begin{figure}[t] \centering
\subfigure[]{\includegraphics[width=3.5cm]{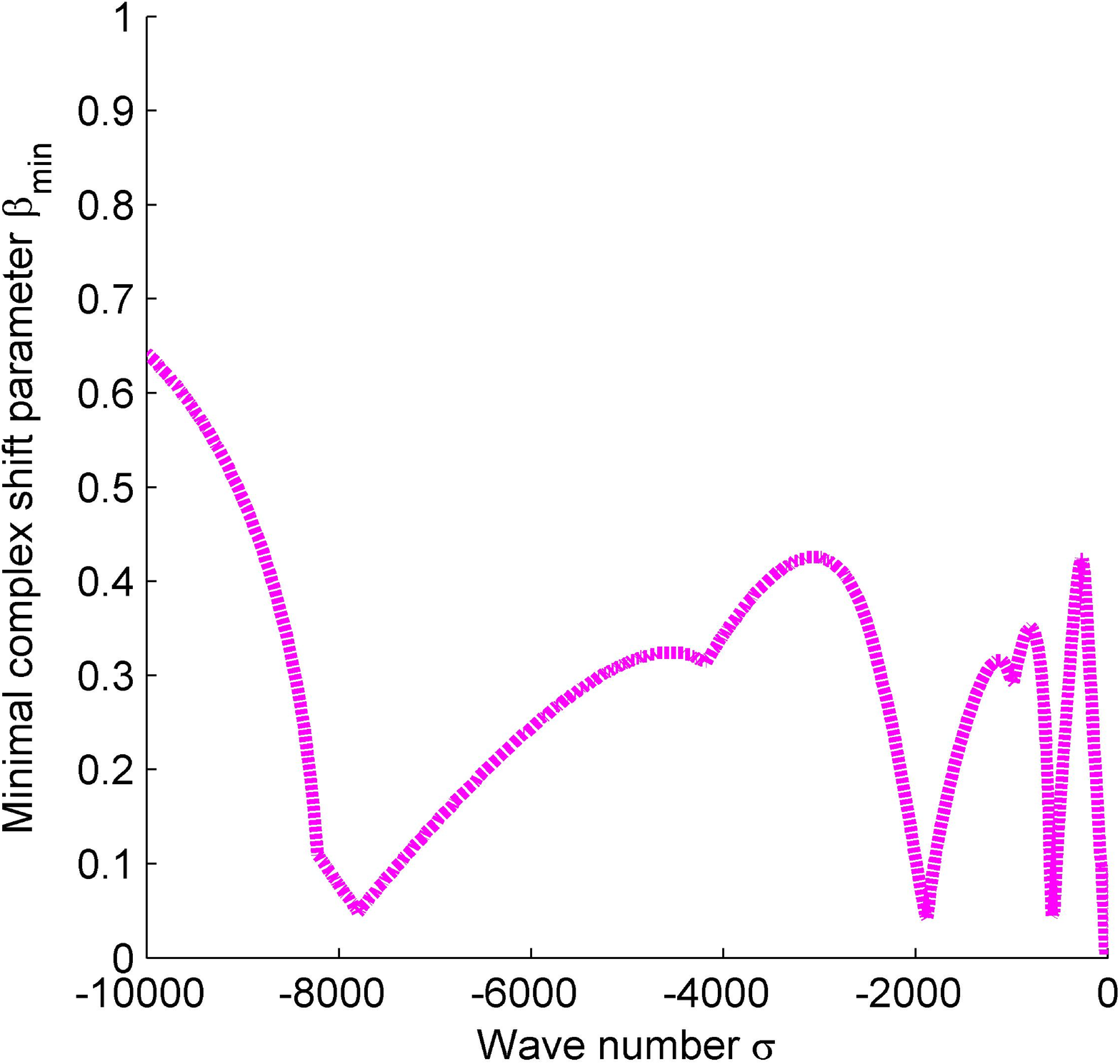}}
\subfigure[]{\includegraphics[width=3.5cm]{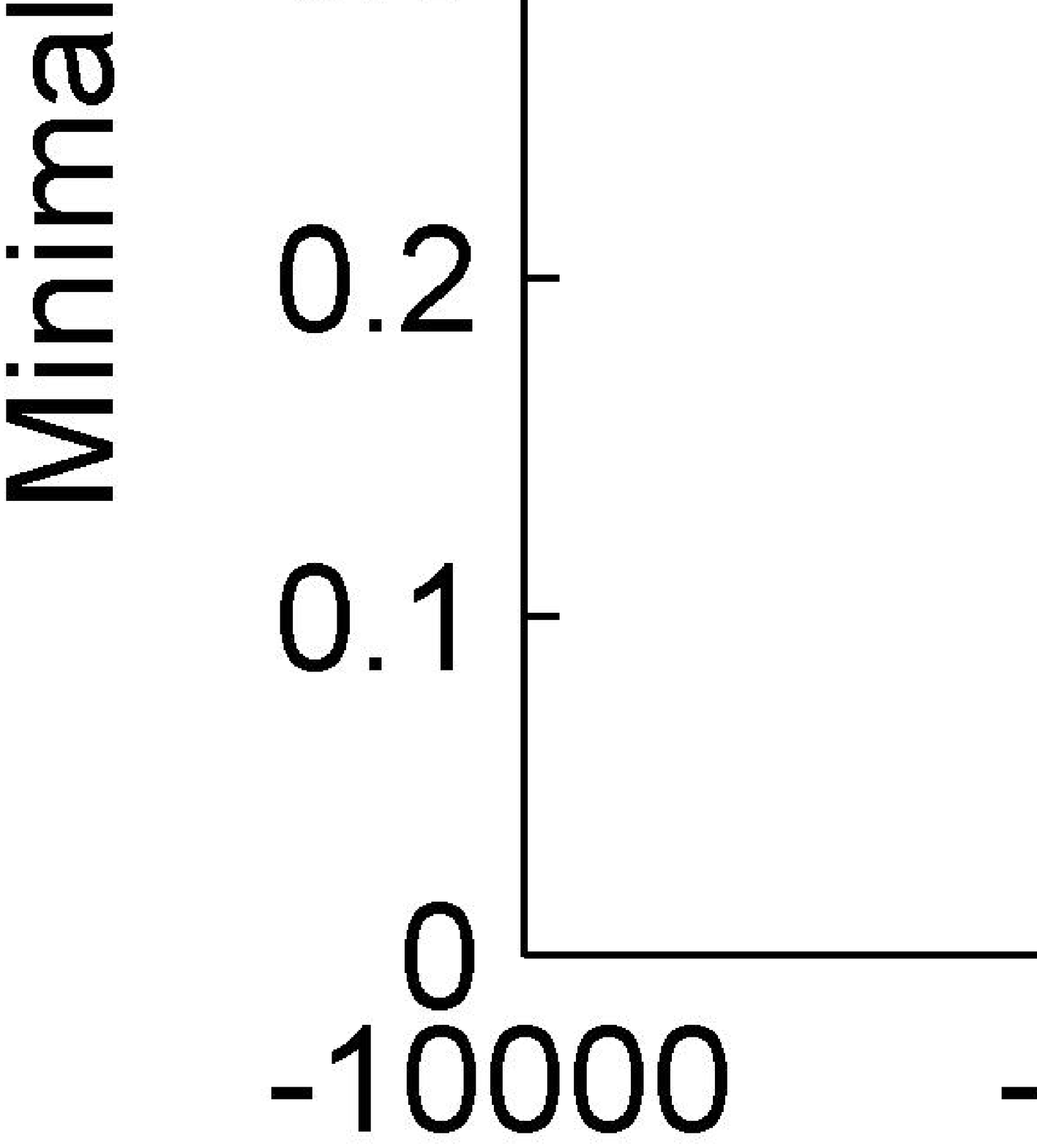}}
\subfigure[]{\includegraphics[width=3.5cm]{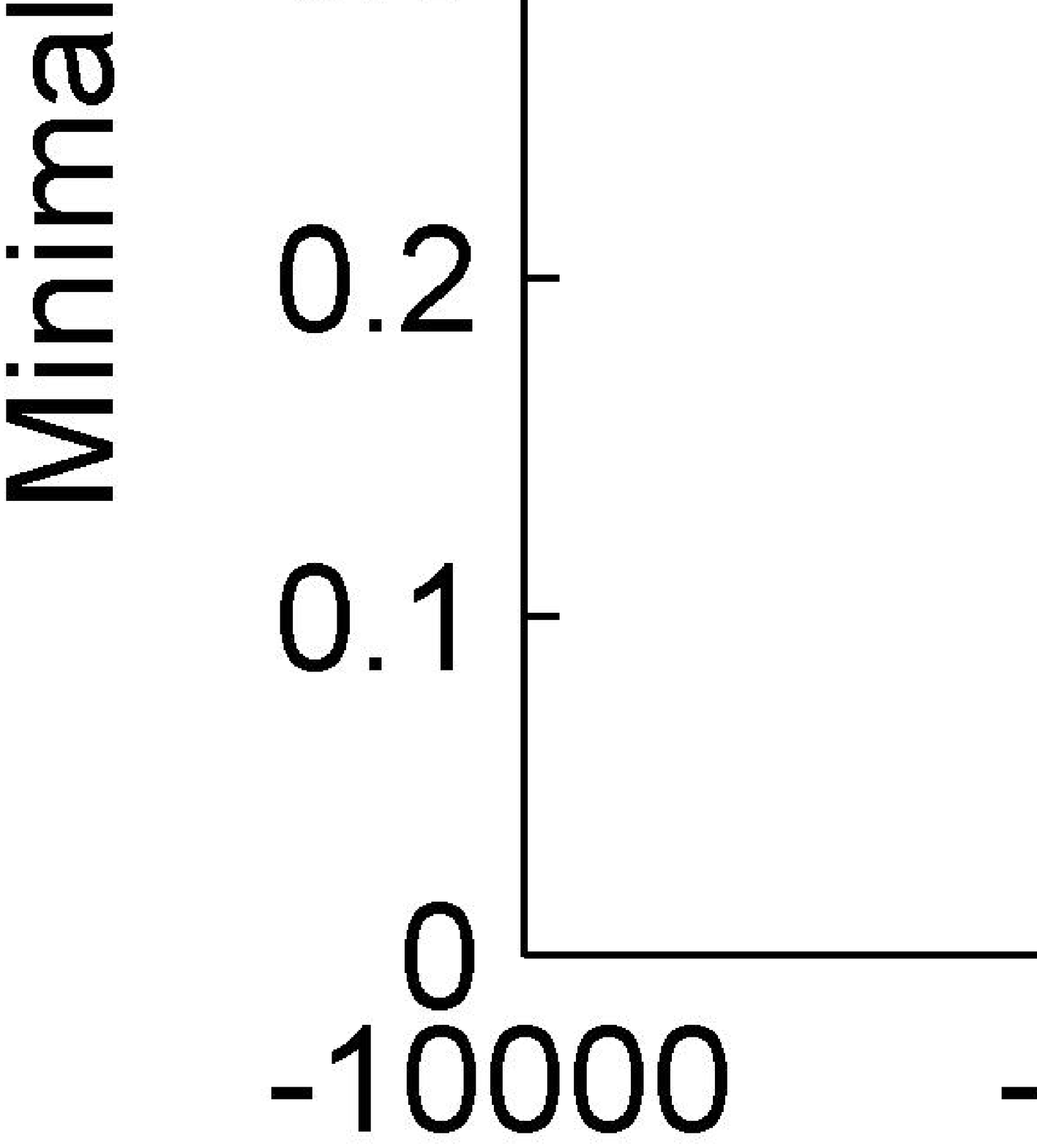}}
\subfigure[]{\includegraphics[width=3.5cm]{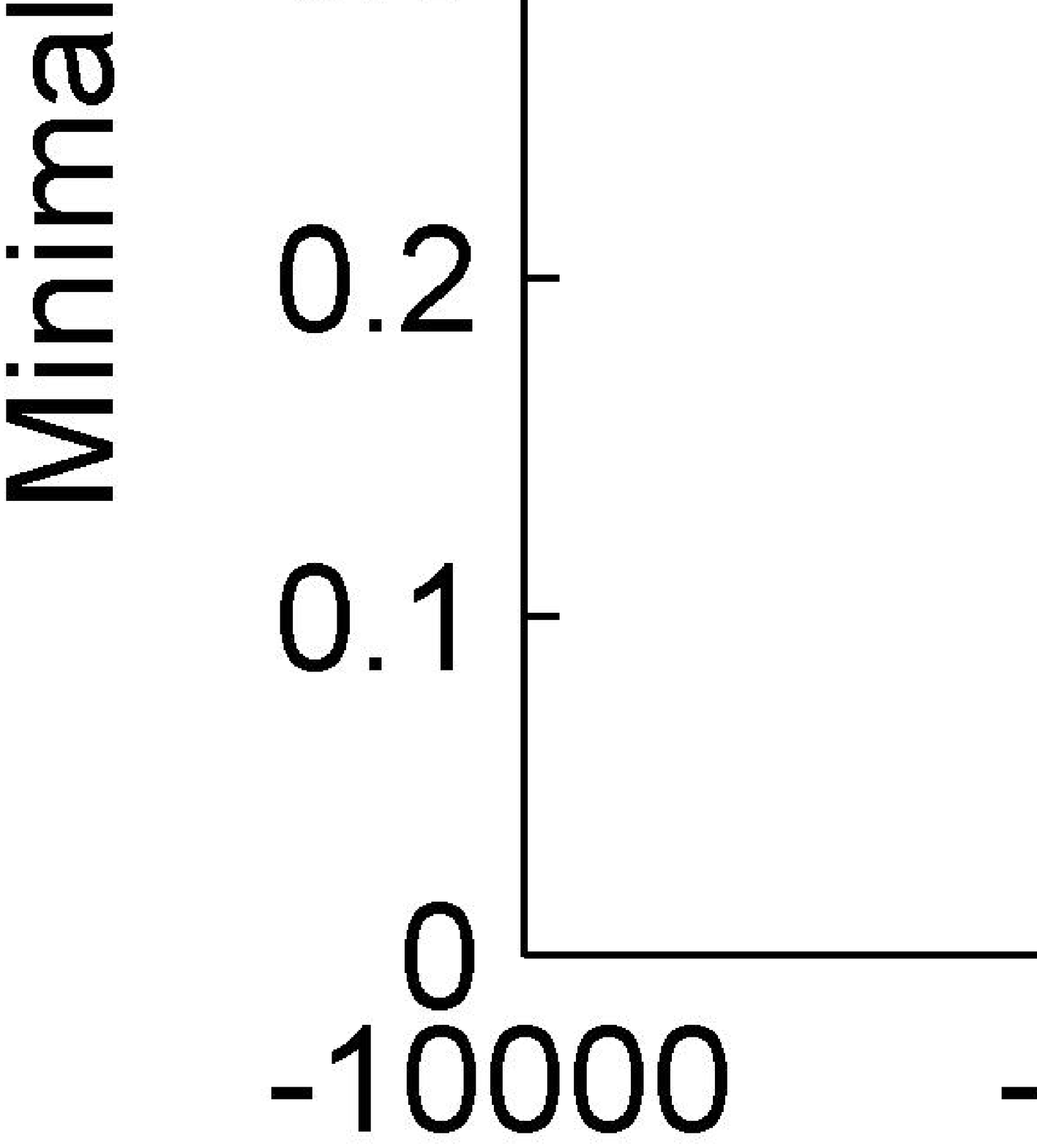}}
\captionsetup{width=14.1cm}
\caption[The minimal complexe shift parameter for different values of $\nu$]{\textsl{The minimal complex shift parameter $\beta_{\min}$ from the 2D 4-grid analysis with $N=64$ and $\omega = 2/3$ as a function of the wavenumber $\sigma$, where $\nu = 1$ (a), $\nu = 2$ (b), $\nu = 3$ (c) and $\nu = 4$ (d).}}
\label{fig:2Dmultismooth}
\end{figure}

Note that all previous results were constructed under the assumption that only one pre- or postsmoothing step is applied, i.e.~$\nu = \nu_1 + \nu_2 = 1$. However, often multiple smoothing steps are used to obtain a more accurate or faster converging iterative solution to the given problem. The minimal complex shift obviously depends on the number of smoothing steps $\nu$, cf. (\ref{eq:propM}). This observation is depicted in Figure \ref{fig:2Dmultismooth}, which shows the $\beta_{\min}$-curve for $\nu = 1,\ldots,4$. The instability of the $\omega$-Jacobi smoother operator, caused by divergence of the smoothest modes as described in \cite{elman2002multigrid} and \cite{brandt1986multigrid}, requires the complex shift to rise significantly for some wavenumbers $\sigma$ when applying multiple smoothing steps. Altering the number of smoothing steps clearly has a significant impact on the $\beta_{\min}$-curve. As a general tendency, one could state that the value of $\beta_{\min}$ rises (at most) linearly as a function of the number of smoothing steps applied. 

\begin{figure}[t] \centering
\subfigure[]{\includegraphics[width=3.5cm]{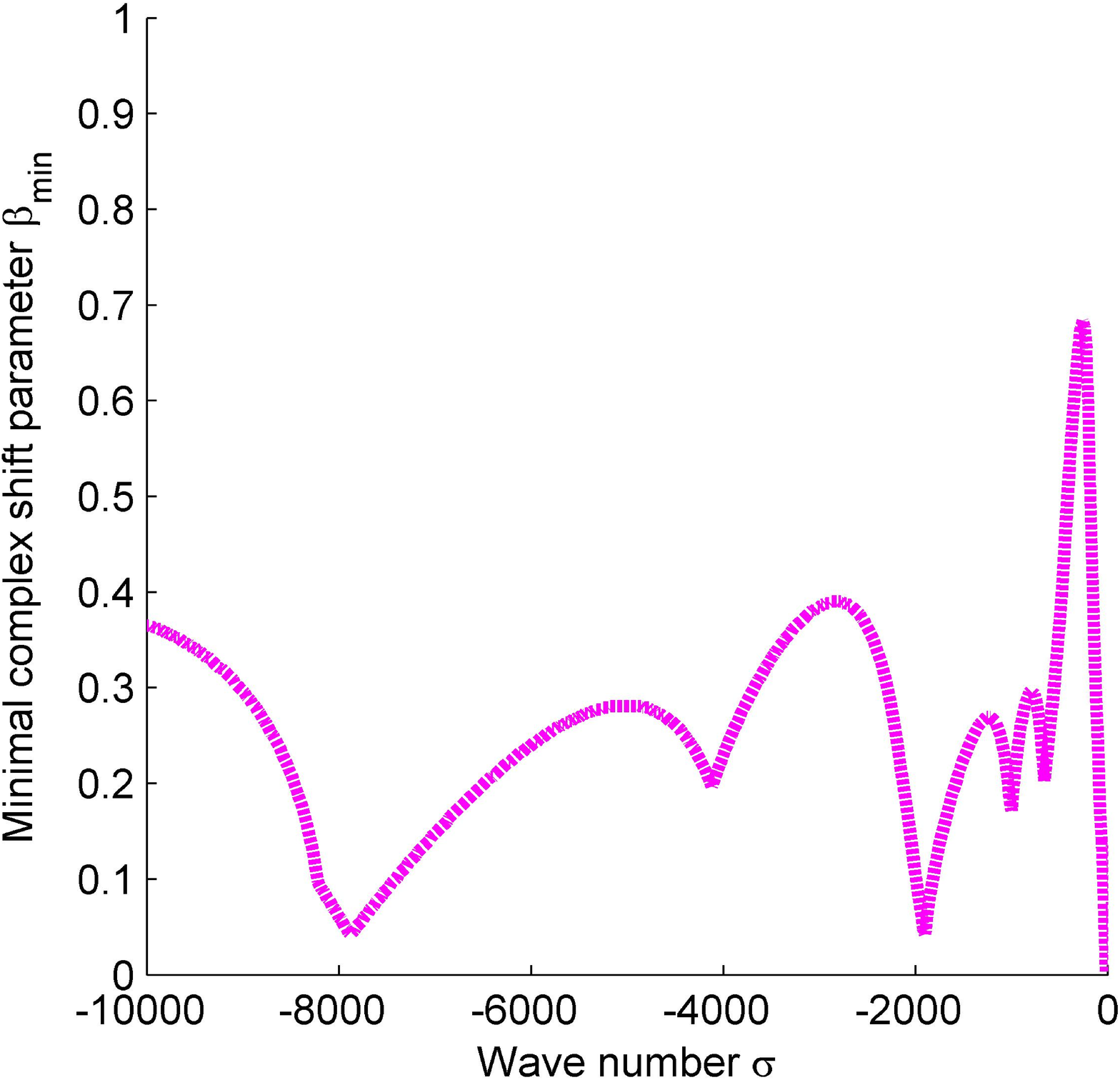}}
\subfigure[]{\includegraphics[width=3.5cm]{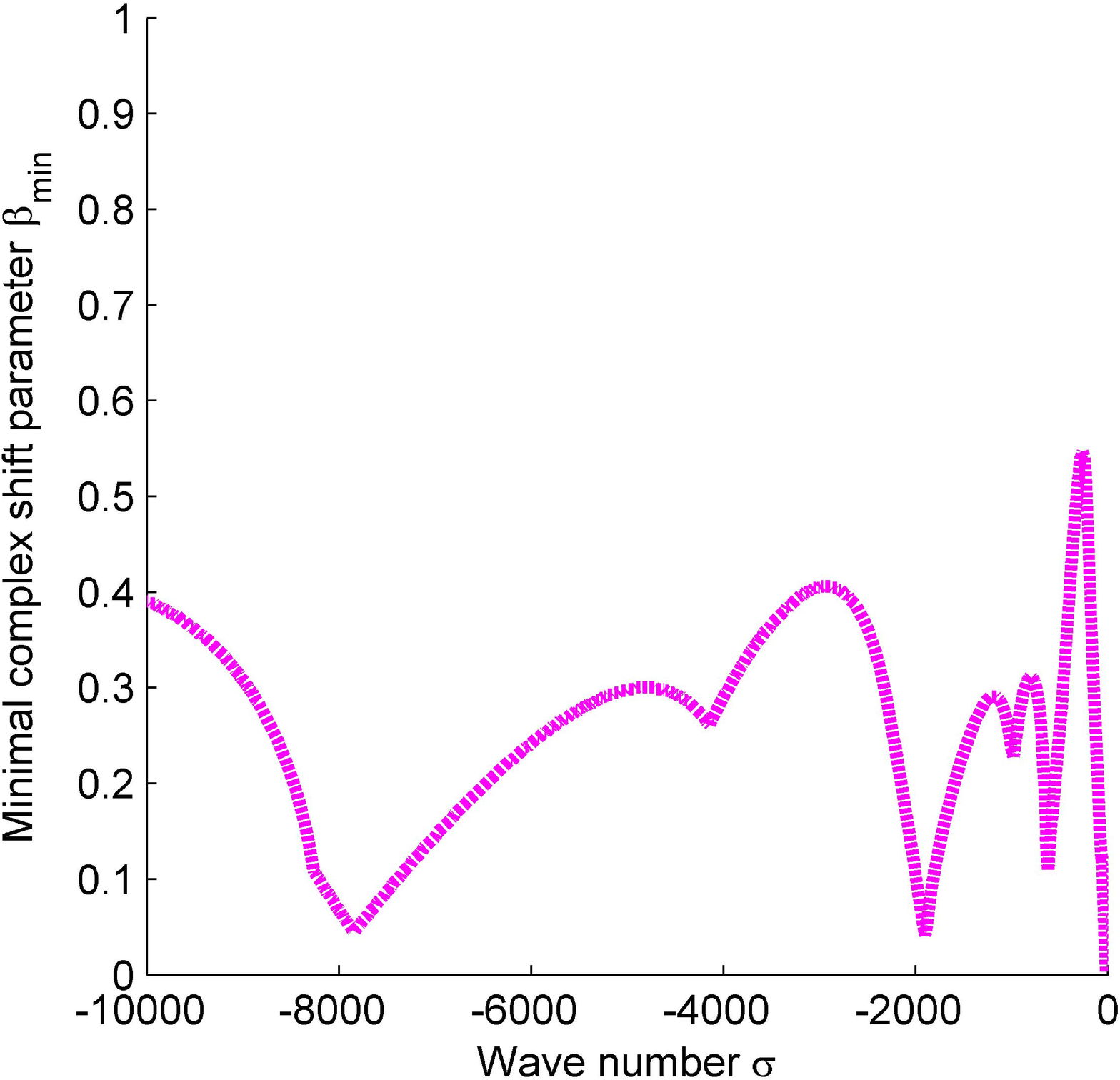}}
\subfigure[]{\includegraphics[width=3.5cm]{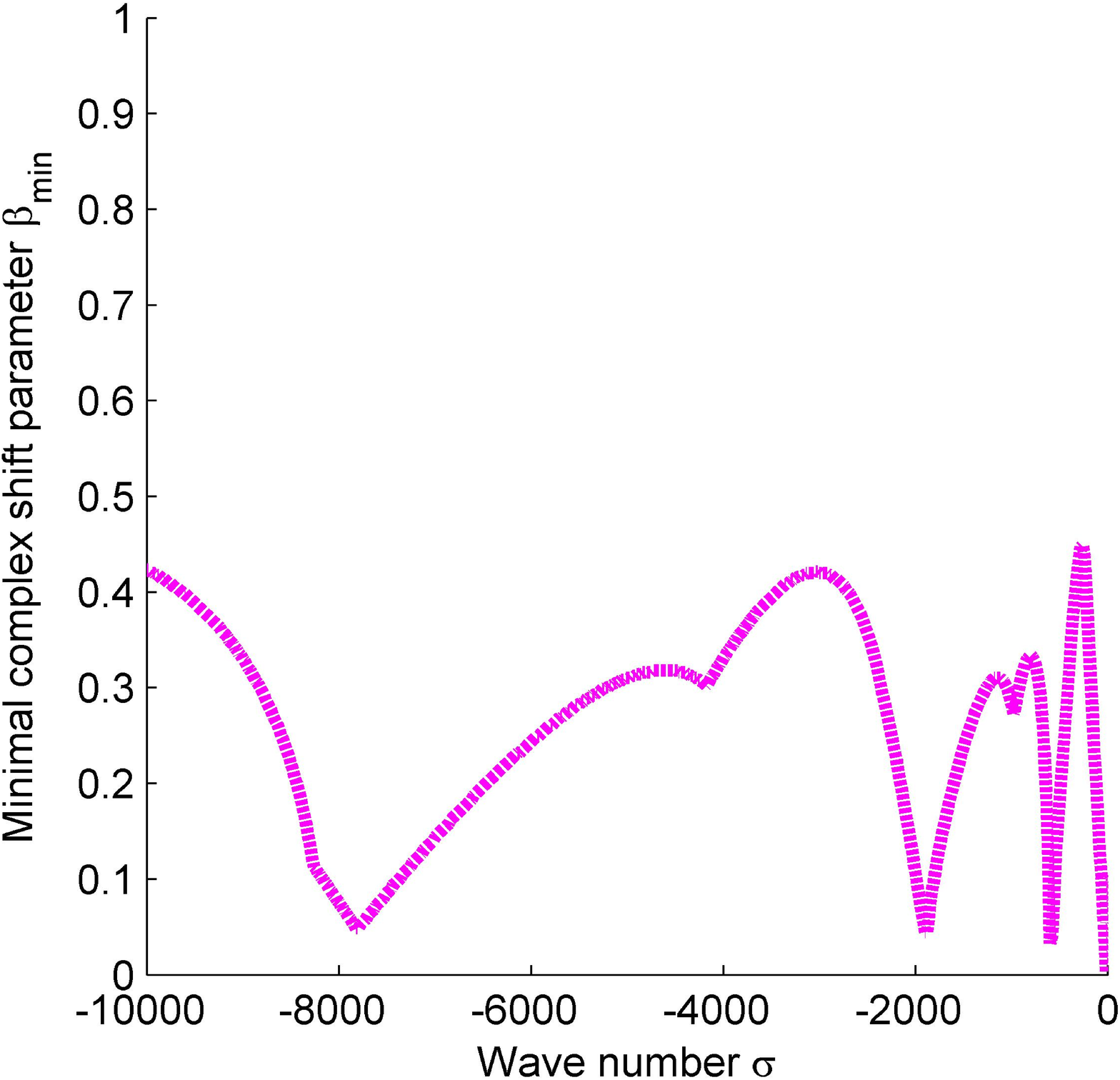}}
\captionsetup{width=14.1cm}
\caption[The minimal complexe shift parameter for different values of $\omega$]{\textsl{The minimal complex shift parameter $\beta_{\min}$ from the 2D 4-grid analysis with $N=64$ and $\nu = 1$ as a function of the wavenumber $\sigma$, where $\omega = 0.2$ (a), $\omega = 0.4$ (b) and $\omega = 0.6$ (c).}}
\label{fig:2Dmultiomega}
\end{figure}

Another parameter of the analysis is the weight $\omega \in [0,1]$ of the Jacobi smoother. When applying the $\omega$-Jacobi smoother, it is convenient to use the standard $\omega = 2/3$ weight, which is known to be optimal for a 1D Poisson problem (see \cite{briggs2000multigrid}). However, for the more general Helmholtz equation with $\sigma \neq 0$, it can be shown (see \cite{elman2002multigrid}) that the optimal Jacobi weight for the 1D problem is given by
\begin{equation}
\omega_{opt} = \frac{2+\sigma h^2}{3+\sigma h^2},
\end{equation}
implying the smoother weight value should be smaller when considering larger values of $|\sigma|$. As shown by Figure \ref{fig:2Dmultiomega}, altering the smoother weight does not have more than a marginal effect on the $\beta_{\min}$-curve, causing it to rise only slightly as $\omega$ increases.\footnote{We remark that function values $\beta_{\min}$ for wavenumbers in the region $|\sigma| < 1000$ should not be taken into account, as the 4-grid scheme is not guaranteed to correctly predict these values - see preceding sections.}

Furthermore, a comment should be made on the discretization-dependency of the theoretical curves. Note that the value of $\beta_{\min}$ is indeed dependend on the finest-mesh stepsize $h = h_1$. However, it is clear from the Fourier symbol calculations that this value never appears separately from the wavenumber $\sigma$, i.e.~$\beta_{\min}$ intrinsically depends on the product $\sigma h^2$. Hence, the value of $\beta_{\min}$ remains unchanged as long as $\sigma h^2$ is constant. In practice, this implies that the $\beta_{\min}$-curve for $N$ evaluated in $\sigma$ has the exact same value as the curve for $2N$ evaluated in $4\sigma$. Doubling the number of finest-grid points therefore implies a stretching of the $\beta_{\min}$-curve over the wavenumber domain by a factor of four. Vice versa, the value of $\beta_{\min}$ in $\sigma/4$ for a $N/2$ discretization is identical to the value of the $\beta_{\min}$-curve in $\sigma$ for $N$ gridpoints, yielding an easy theoretical lower limit for the value of $\beta$ when considering a multigrid V-cycle on a coarser level.

Although the analysis in this paper is restricted to the constant-$k$ model problem given by (\ref{eq:mod1}), the results displayed here can be easily extended to cover more realistic space-dependent $k$ settings. Indeed, for problems with a moderately space-dependent wavenumber like e.g.~the regionally space-dependent \emph{Wedge model}, see \cite{erlangga2006novel}, the previous analysis may be conducted on each distinct stencil, resulting in each domain region being associated with a different minimal complex shift $\beta_{\min}$. Combining these lower limits, an indisputably safe choice for the complex
shift parameter $\beta$ would be the largest possible regional $\beta_{\min}$ appearing in the problem, which leads to a stable multigrid scheme in every region of the domain.

To conclude this section, we present a short discussion on the choice of $\beta_{\min}$ that is currently being used in the literature. A widely accepted choice for the complex shift parameter is $\beta = 0.5$, which was first introduced in \cite{erlangga2006novel} and has been used ever since by many researchers in the field. In this paper by Erlangga, Oosterlee and Vuik, one reads that \\

\begin{minipage}{13.5cm}
\emph{``The preconditioner of choice in this paper is based on the parameters $(\beta_1, \beta_2) = (1, 0.5)$. (\ldots) For values $\beta_2 < 0.5$ it is very difficult to define a satisfactory converging multigrid F(1,1)-cycle with the components at hand. They are therefore not considered. (\ldots) From the results in Table 6 we conclude that the preferred methods among the choices are the preconditioners with $\beta_1 = 1$. (\ldots) Fastest [multigrid preconditioned Krylov] convergence is obtained for $(\beta_1,\beta_2) = (1,0.5)$.''}
\end{minipage}\\ \vspace{2mm}

\noindent where $\beta_1$ and $\beta_2$ refer to the real and complex shift parameters, designated in this text by $\alpha$ and $\beta$ respectively. Conclusions were drawn from a variety of numerical experiments with a fixed wavenumber $\sigma$ and mesh width satisfying $k h = 0.625$. The results can be compared to the theoretical value for $\beta_{\min}$ suggested by Figure \ref{fig:2Dmultismooth}(b) at $\sigma = -1600$. Note that we apply a V-cycle and use the smoother weight $\omega = 2/3$, contrary to the F-cycle and $\omega = 0.5$ used in \cite{erlangga2006novel}, however these differences only minorly affect the $\beta_{\min}$ curve. It is clear that for all wavenumbers $\sigma \geq -1600$ (indicating the use of at least 10 gridpoints per wavelength), the corresponding minimal complex shift is indeed smaller than 0.5. In the case where $\sigma = -1600$, $\beta_{\min}$ equals $0.2$, suggesting an even smaller complex shift may be used. Hence, respecting the $kh \leq 0.625$ criterion, the choice for $\beta = 0.5$ always guarantees V(1,1)-cycle convergence. The supposed near-optimality of $\beta_{\min}$ with respect to the number of Krylov iterations is discussed in the next section.

\newpage

\section{Numerical results}

In this section, we present some experiments that will be used to confirm and validate the theoretical results obtained in Section 2, as well as provide the reader with some valuable insight on the definition choice of the minimal complex shift as stated above. The value of $\beta_{\min}$ predicted by the Local Fourier Analysis curves will indeed prove a valuable lower limit for $\beta$ in practical applications. Note that some of the experimental figures obtained in this section are mainly intended as a validation of the theoretical results, and their further use in practical applications is rather limited.

\begin{figure}[t] \centering
\subfigure[]{\includegraphics[width=4.5cm]{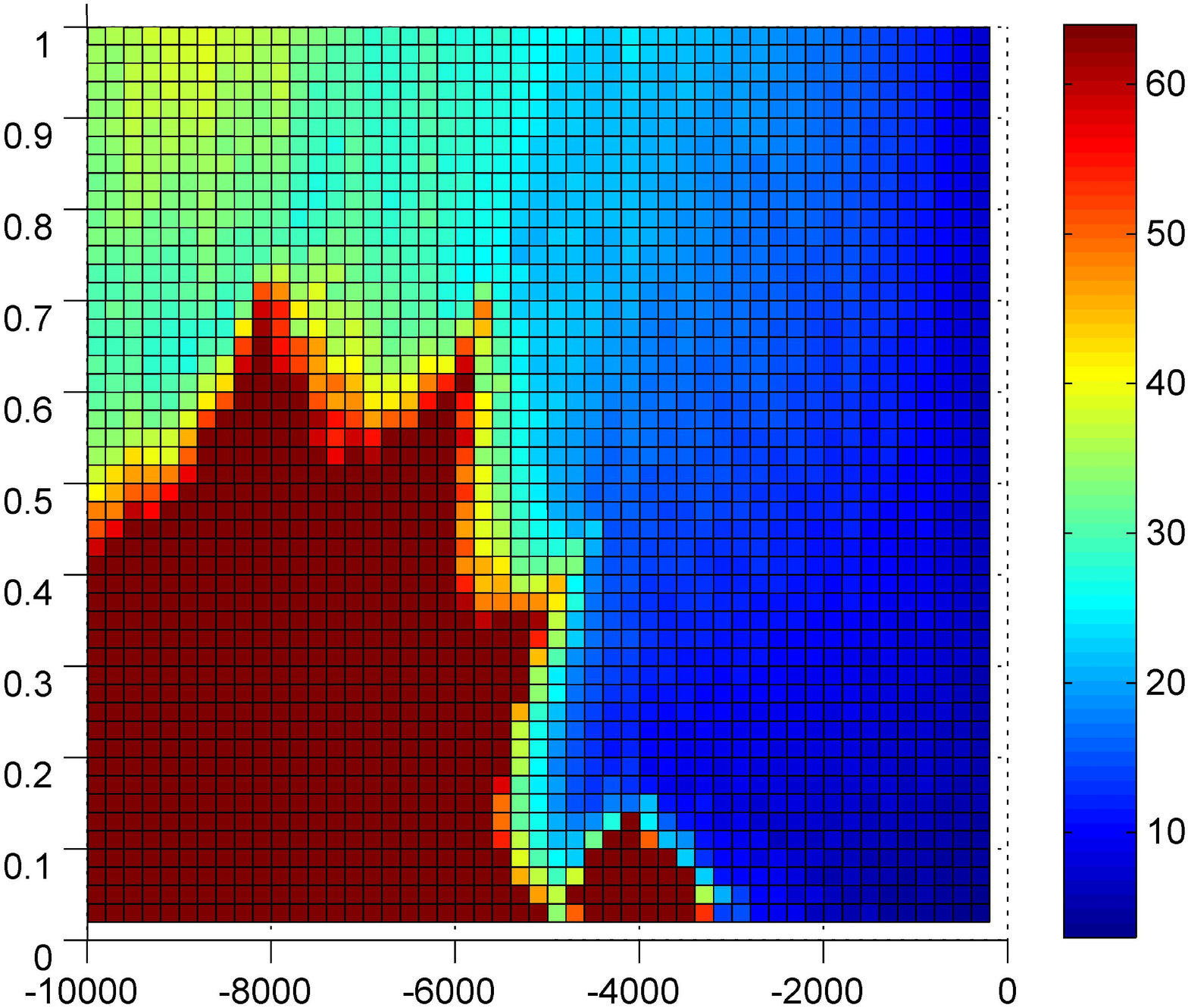}}
\subfigure[]{\includegraphics[width=4.5cm]{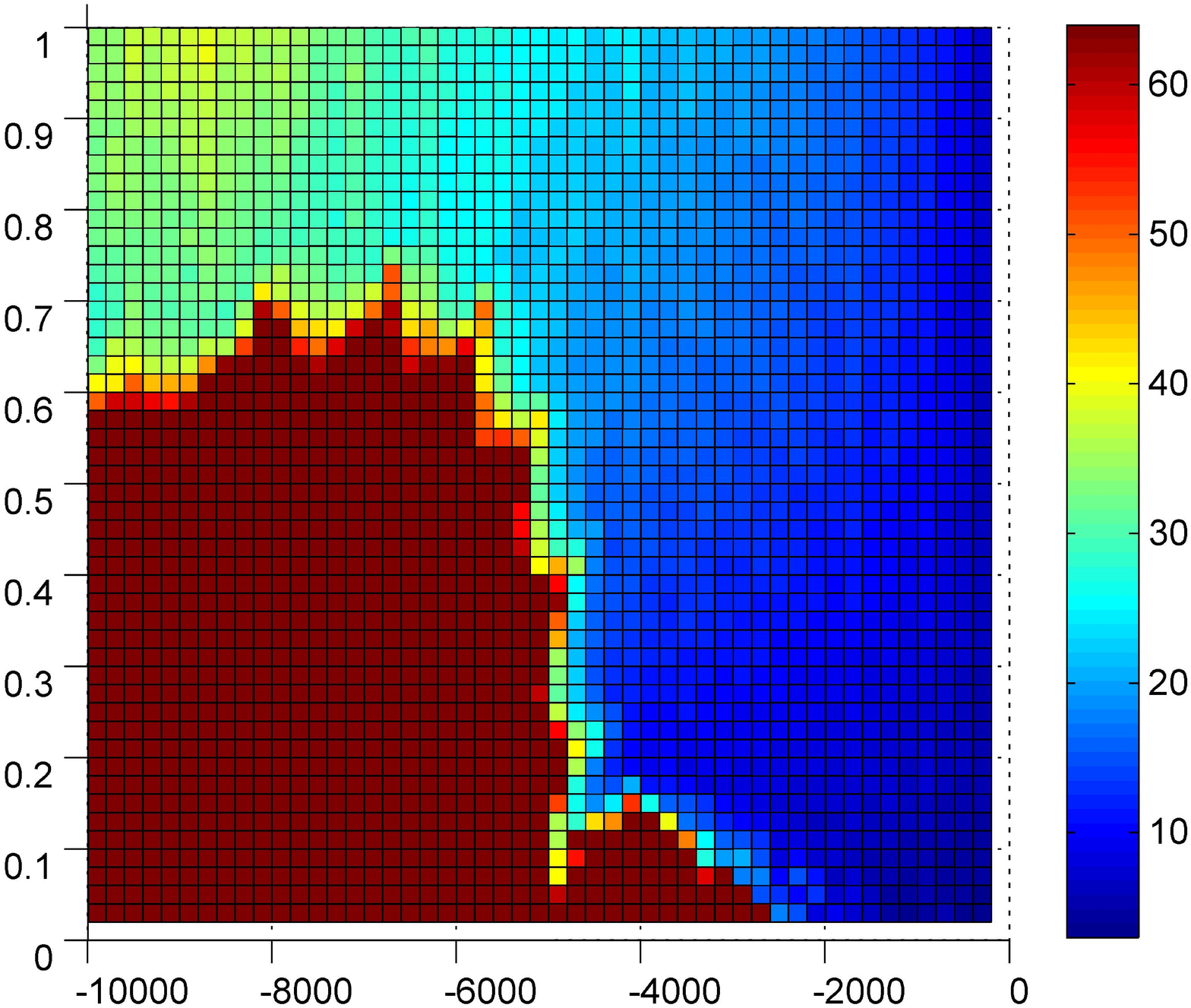}}
\subfigure[]{\includegraphics[width=4.5cm]{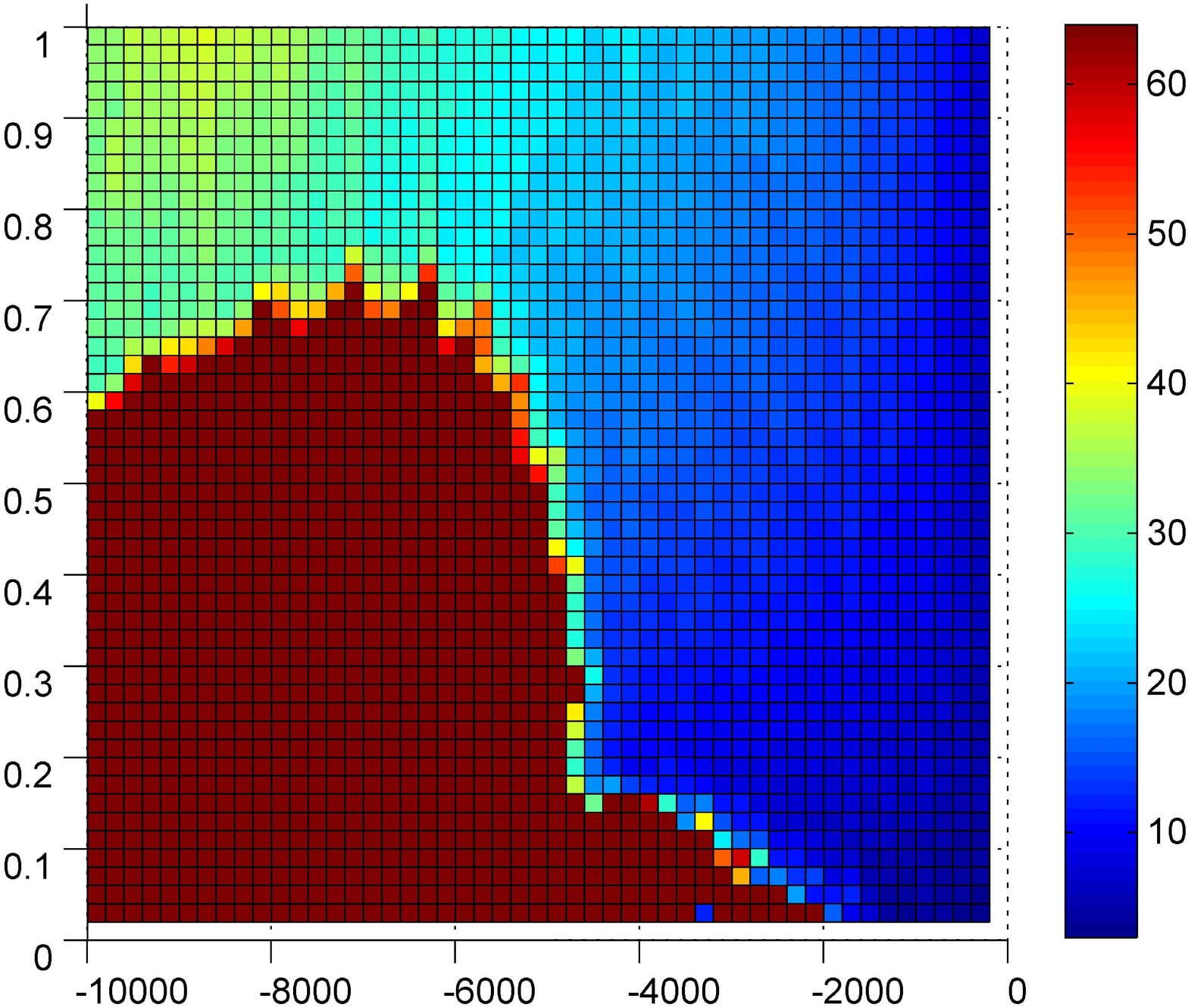}}
\captionsetup{width=14.1cm}
\caption[1D two-grid preconditioned BiCGStab iteration count]{\textsl{Two-grid-preconditioned Krylov method iteration count (colour) as a function of the wavenumber $\sigma$ and the complex shift $\beta$ for the 1D model problem with $N=64$. The applied method is BiCGStab with $\mu = 10$ (a), $\mu = 20$ (b) and $\mu = 30$ (c).}}
\label{fig:1DBiCGTG}
\end{figure}

\begin{figure}[t] \centering
\subfigure[]{\includegraphics[width=7cm]{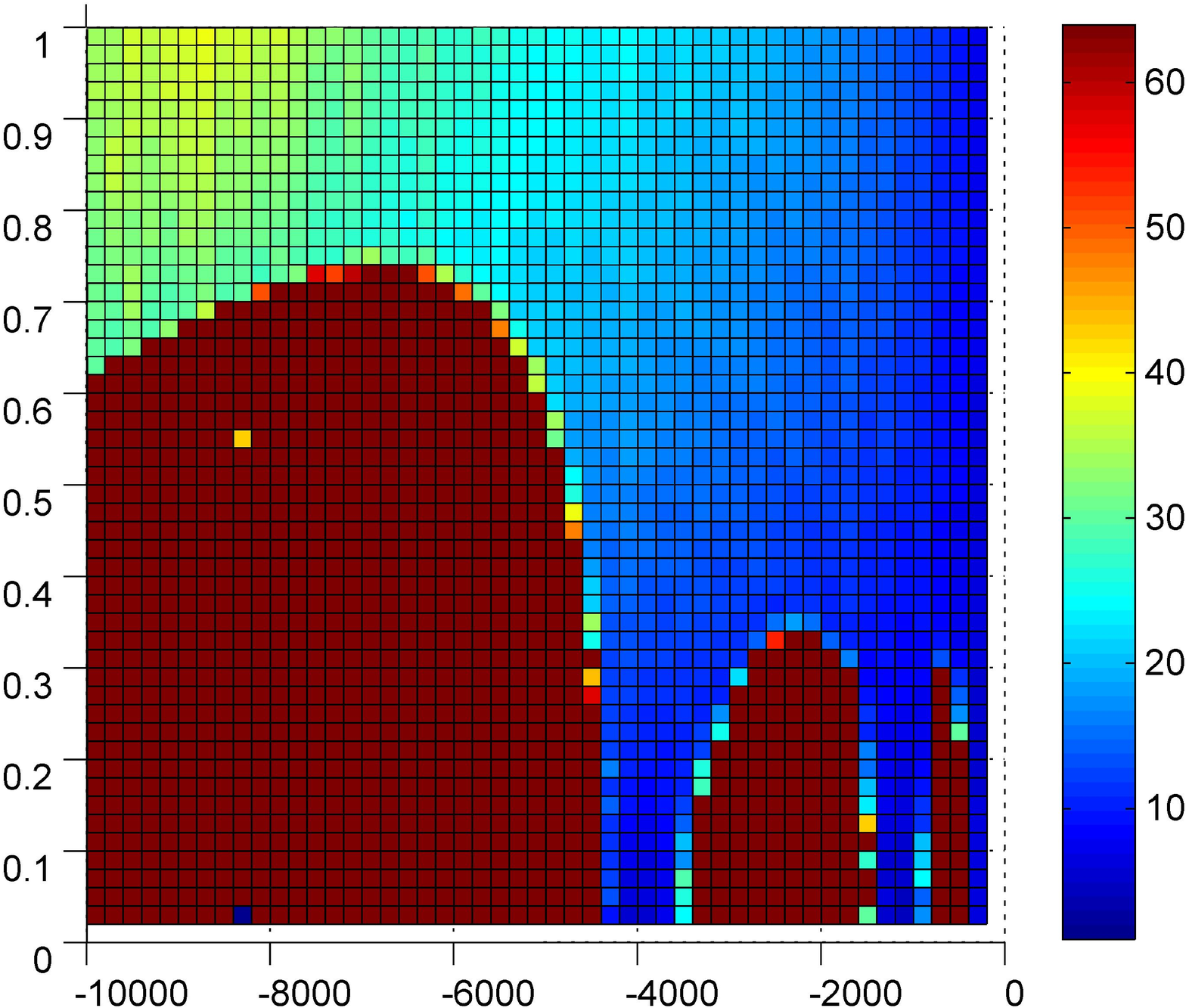}}
\subfigure[]{\includegraphics[width=7cm]{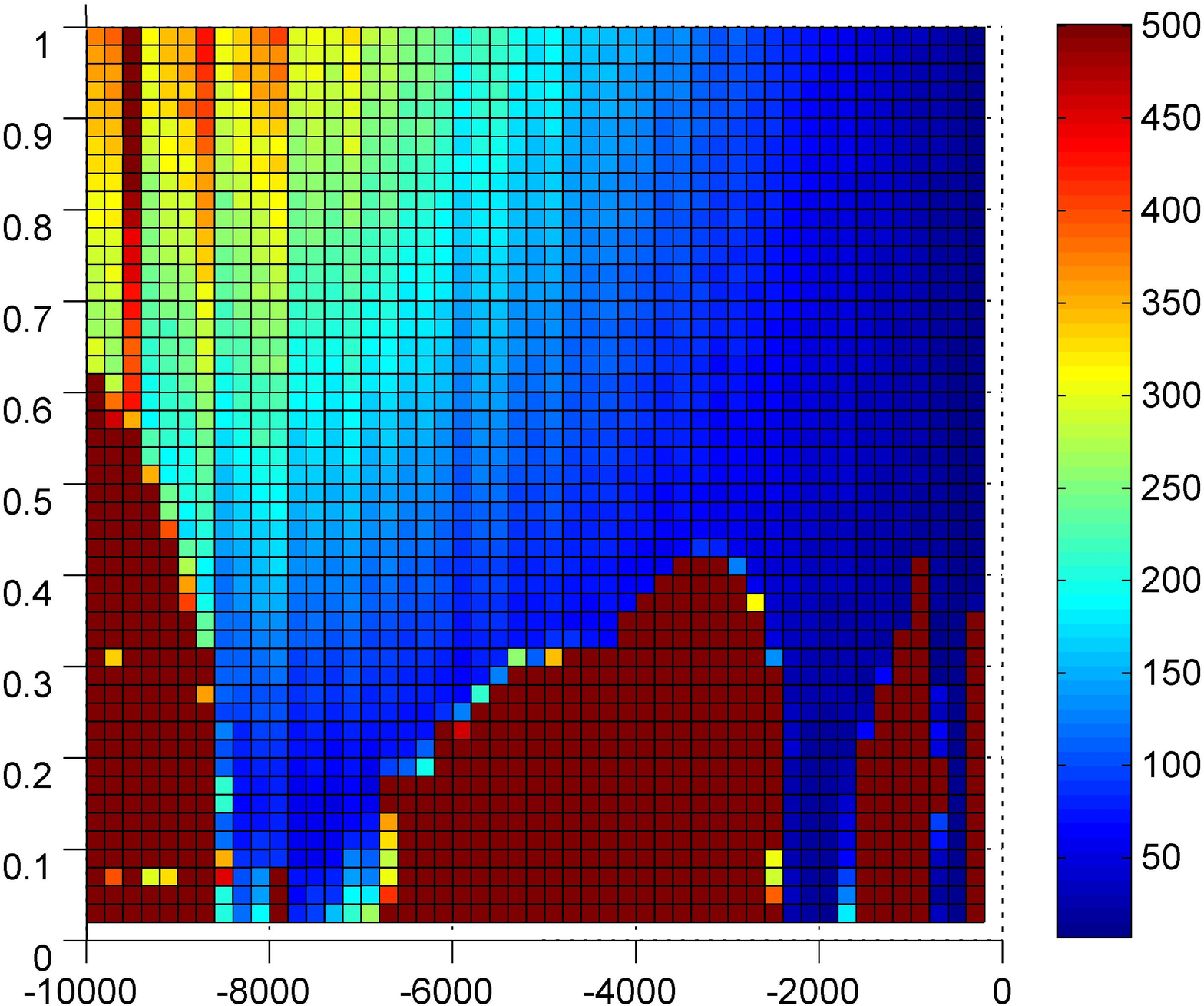}}
\captionsetup{width=14.1cm}
\caption[1D and 2D BiCG V-cycle preconditioned BiCGStab iteration count]{\textsl{V(1,0)-cycle-preconditioned Krylov method iteration count (colour) as a function of the wavenumber $\sigma$ and the complex shift $\beta$ for the 1D (a) and 2D (b) model problems, both with $N=64$. The applied method is BiCGStab with $\mu = 100$ (a) and $\mu = 60$ (b).}}
\label{fig:1D2DBiCGV}
\end{figure}

\begin{figure}[t] \centering
\subfigure[]{\includegraphics[width=6cm]{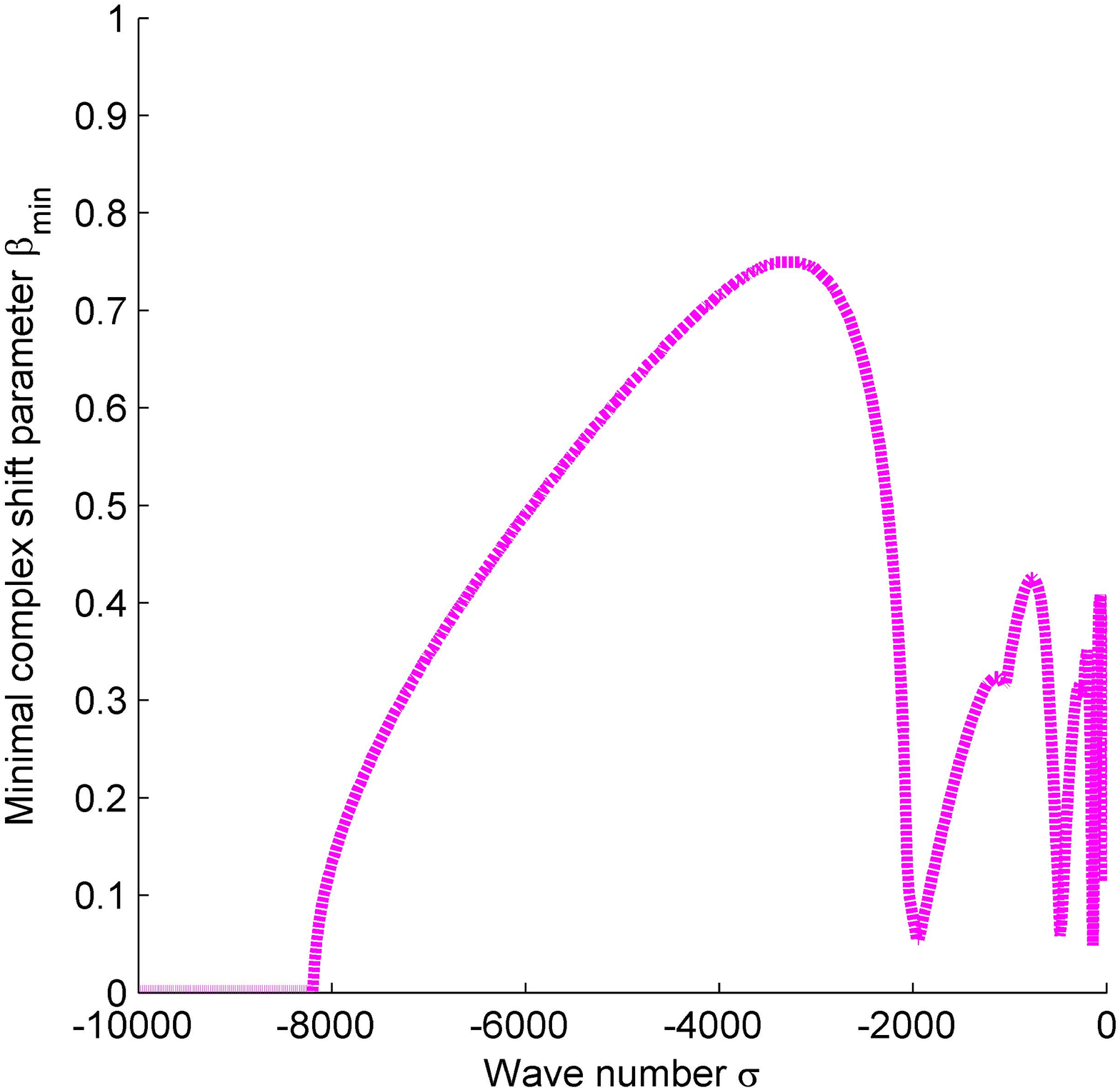}}
\subfigure[]{\includegraphics[width=7.4cm]{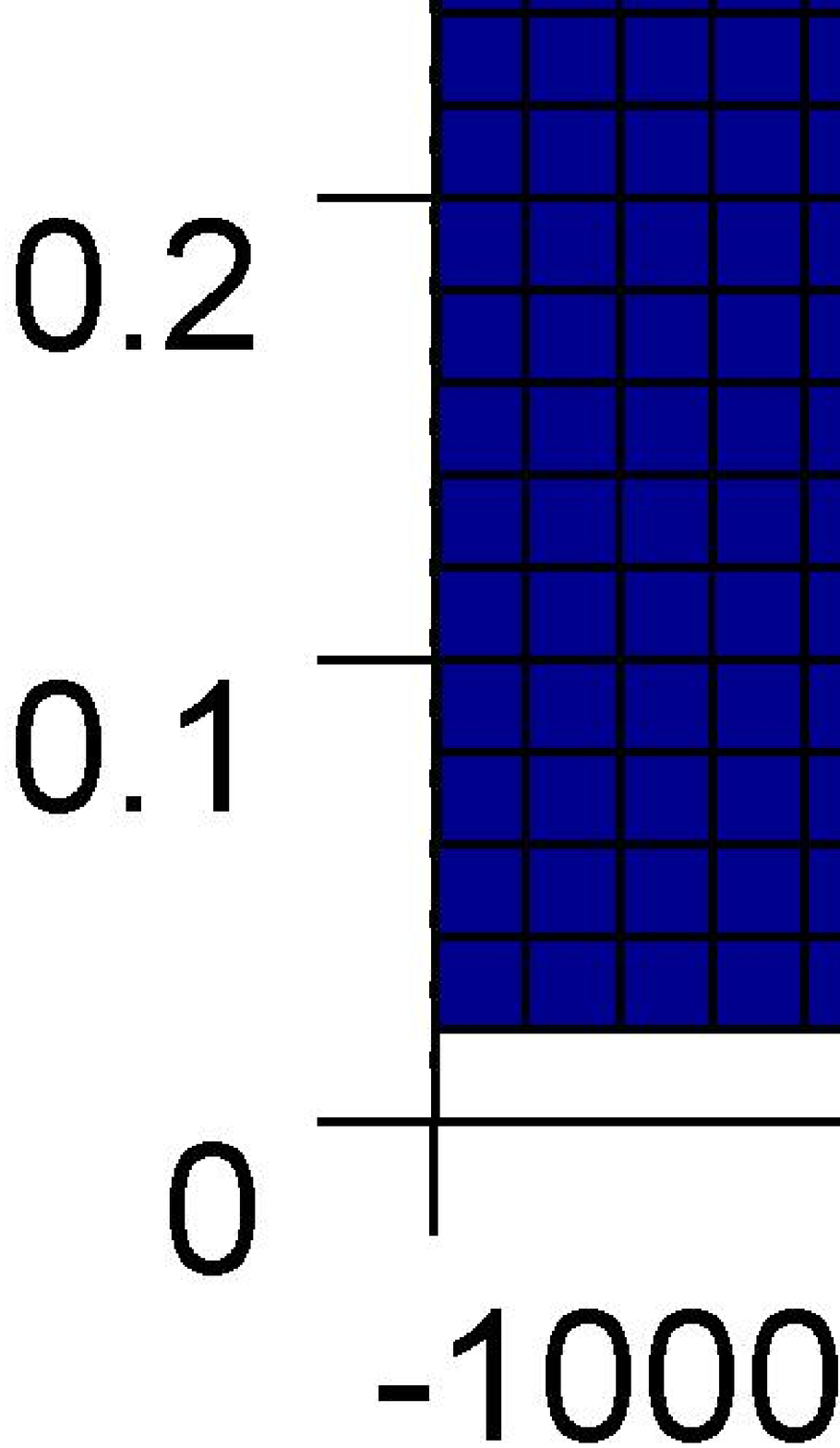}}
\captionsetup{width=14.1cm}
\caption[]{\textsl{V(1,0)-cycle-preconditioned GMRES iteration count with $\mu=30$ (colour) as a function of the wavenumber $\sigma$ and the complex shift $\beta$ for the 2D model problem with $N=32$ (b), in comparison to the theoretical 4-grid LFA complex shift parameter $\beta_{\min}$ (a).}}
\label{fig:2DbetaN32}
\end{figure}

\begin{figure}[t] \centering
\subfigure[]{\includegraphics[width=6cm]{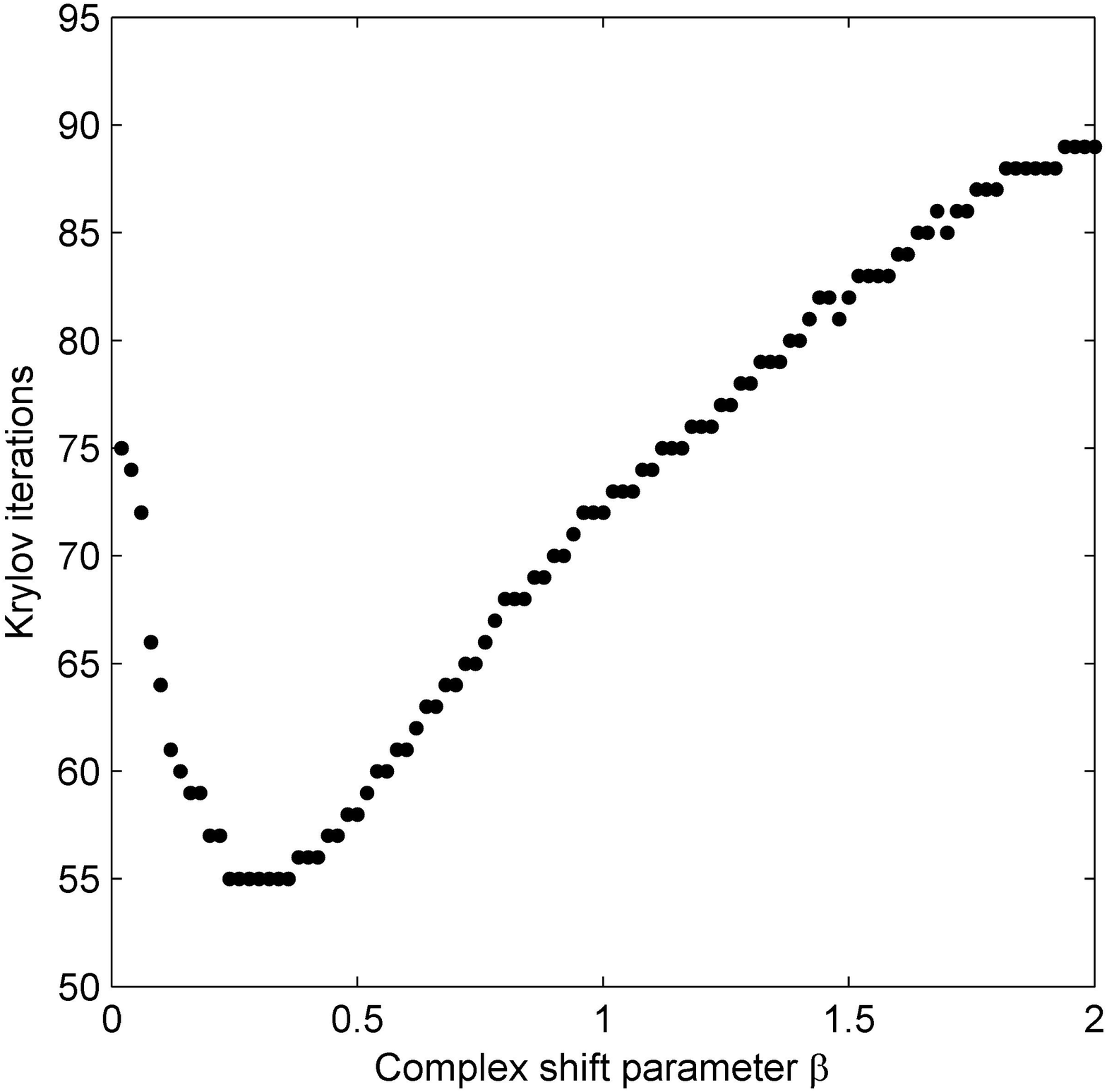}}
\subfigure[]{\includegraphics[width=6cm]{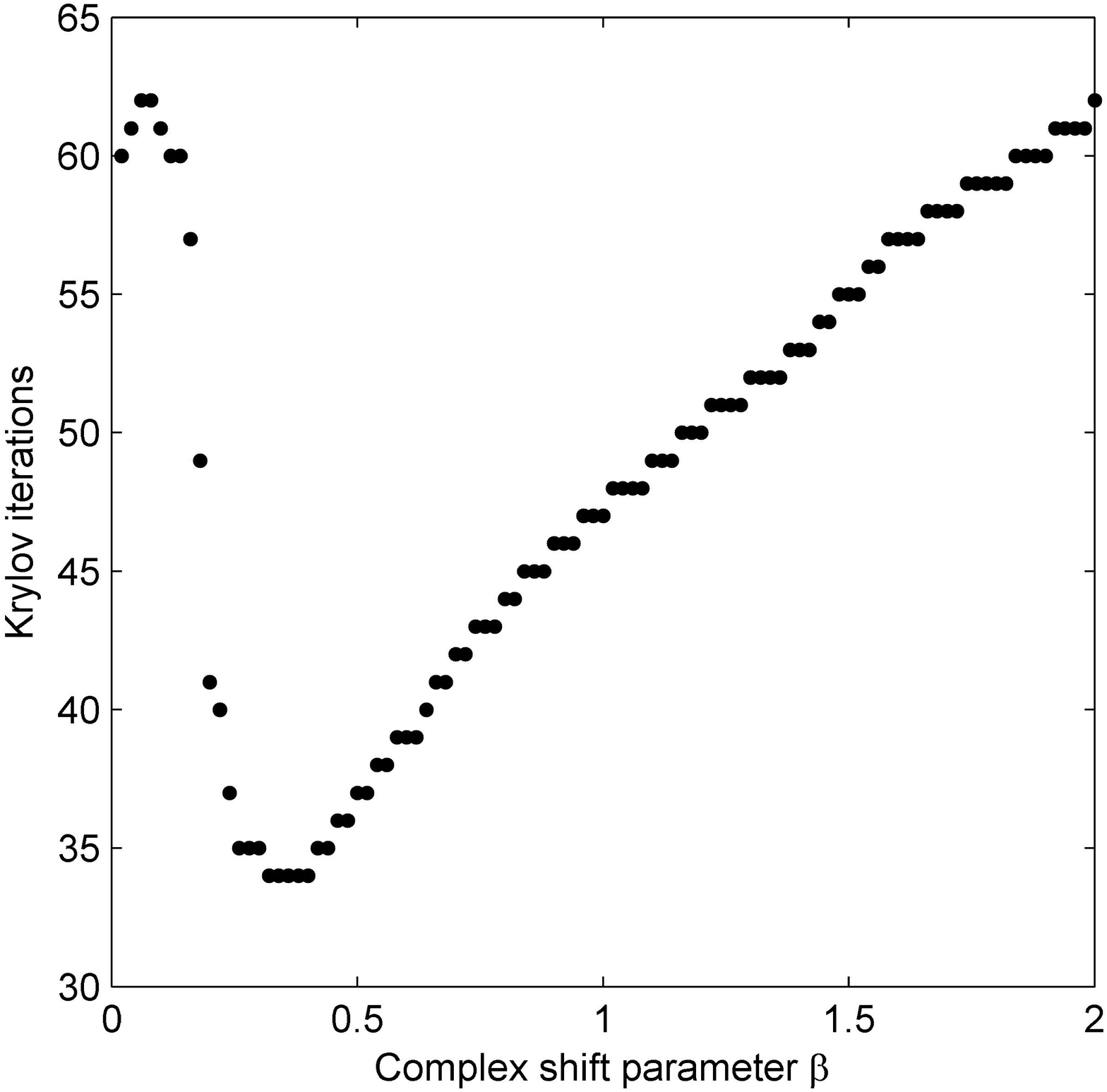}} 
\subfigure[]{\includegraphics[width=6cm]{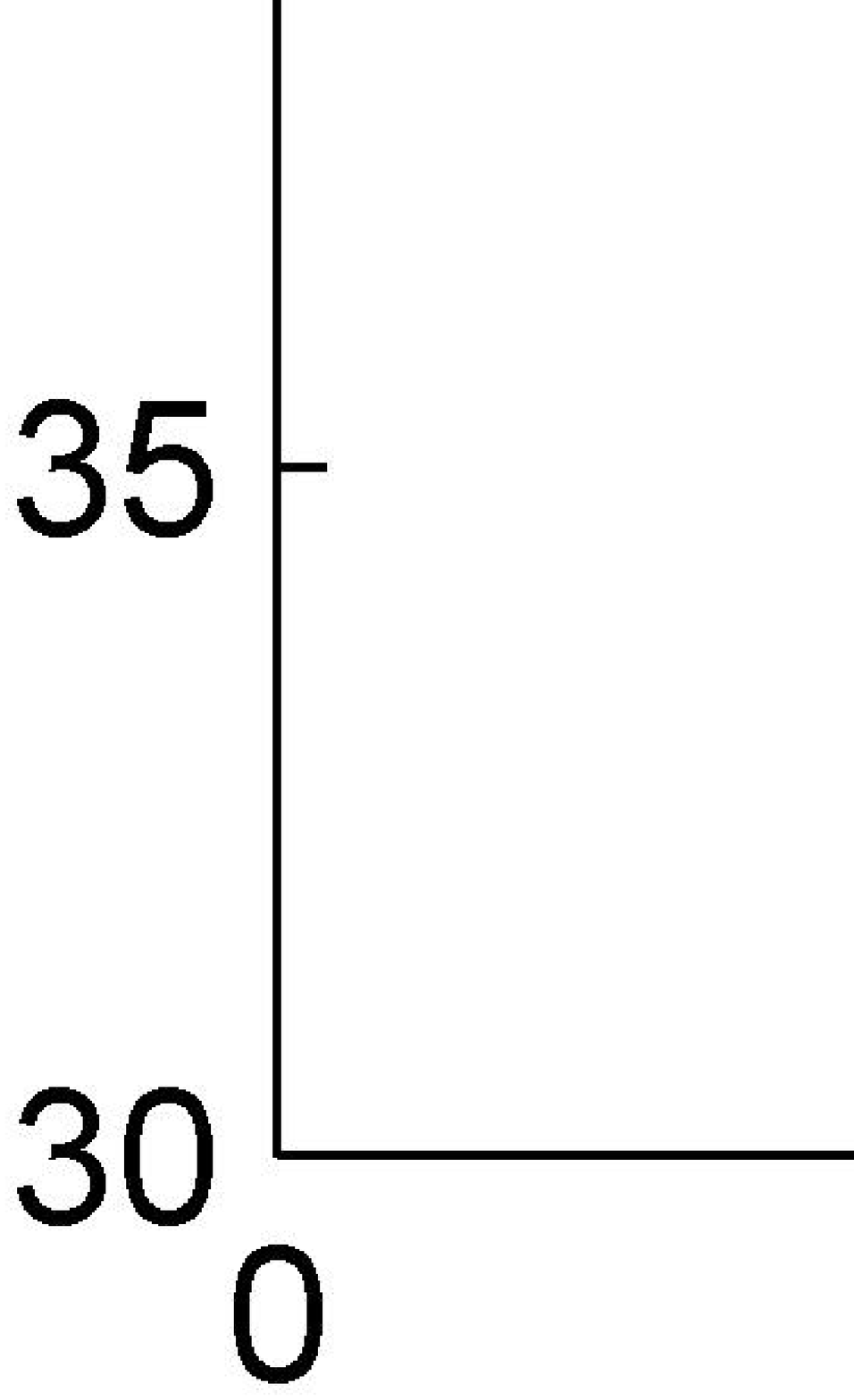}}
\subfigure[]{\includegraphics[width=6cm]{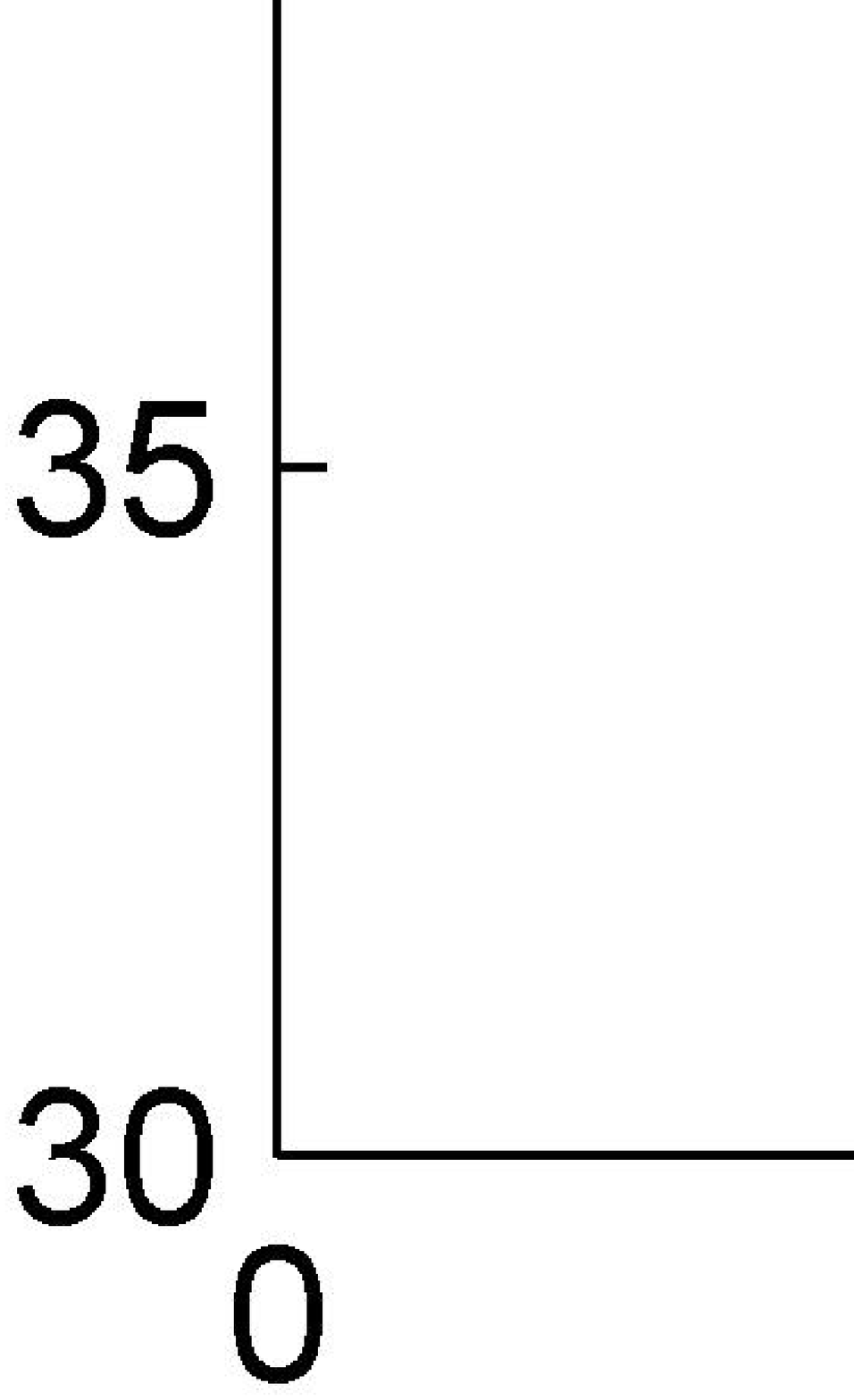}}
\caption[]{\textsl{Two-grid-preconditioned Krylov method iteration
    count for $\sigma = -1000$ in function of the complex shift $\beta$ for the 2D model problem
    with $N=32$. The applied method is GMRES with $\mu=1$ (a), $\mu=3$ (b),
    $\mu=5$ (c) and $\mu=10$ (d). Corresponding minima can be found at $\beta = 0.30$ (a), $\beta = 0.36$ (b), $\beta = 0.40$ (c) and $\beta = 0.42$ (d). \vspace{-1.2cm}}}
\label{fig:2DGMRESTGIters}
\end{figure}

\begin{figure}[t] \centering
\includegraphics[width=7cm]{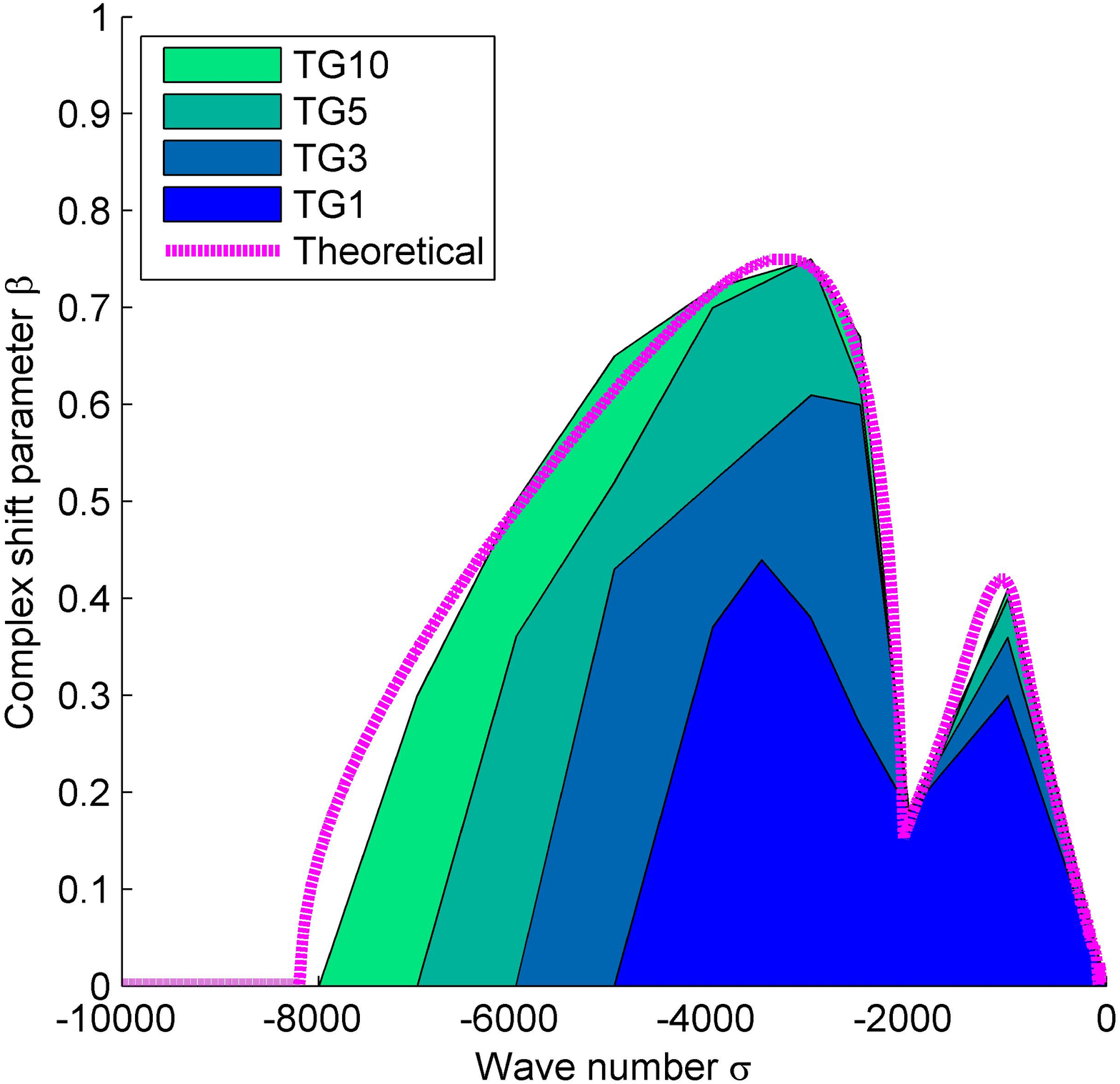}
\captionsetup{width=14.1cm}
\caption[2D two-grid preconditioned GMRES lowest iteration number]{\textsl{Two-grid-preconditioned Krylov method iteration-minimum-$\beta$ as a function of the wavenumber $\sigma$ for the 2D model problem with $N=32$. The applied method is GMRES with $\mu = 1$, $\mu=3$, $\mu=5$ and $\mu=10$.}}
\label{fig:2DGMRESTGItersLowest}
\end{figure}

\begin{figure}[t] \centering
\subfigure[]{\includegraphics[width=4.5cm]{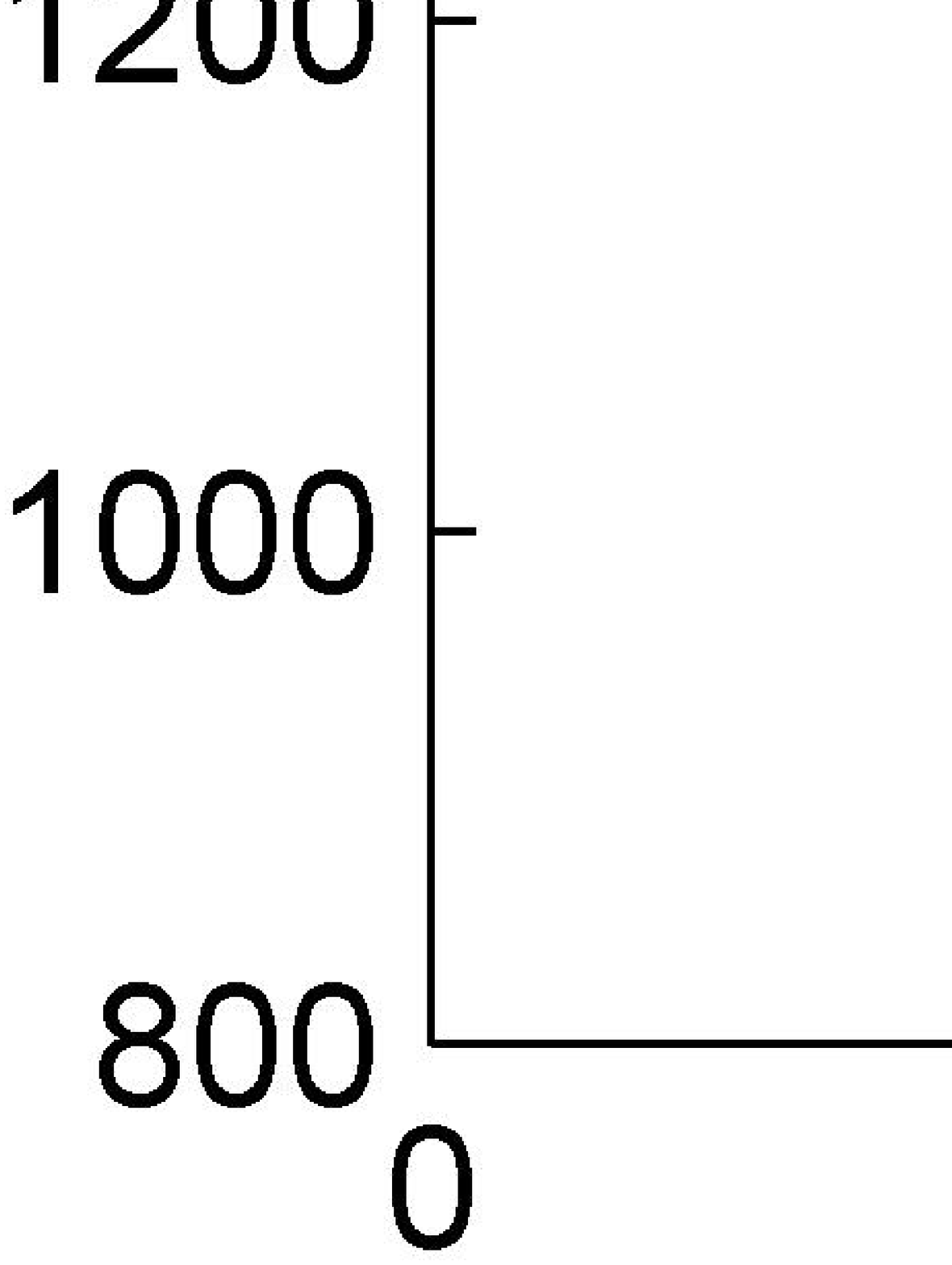}}
\subfigure[]{\includegraphics[width=4.5cm]{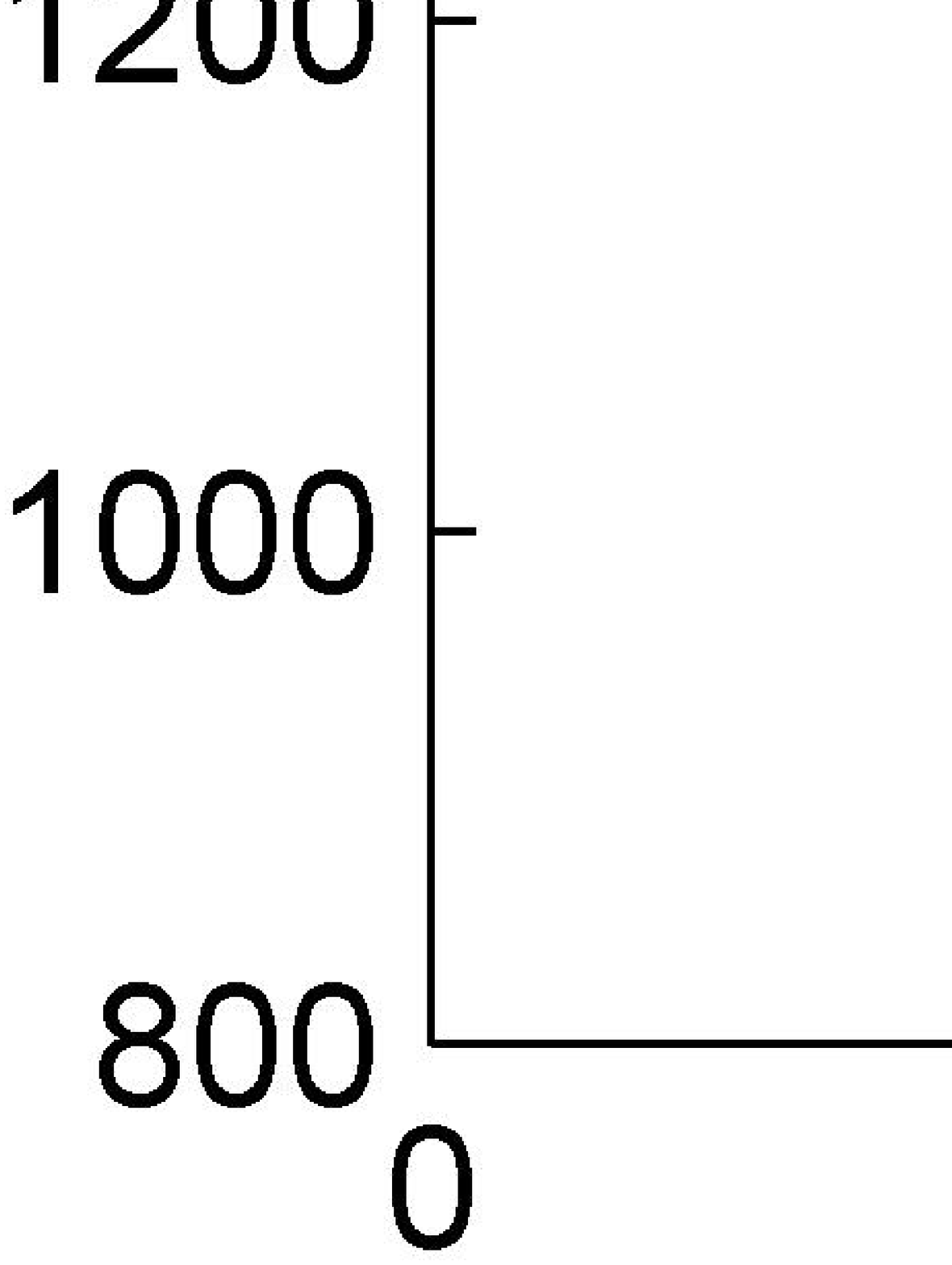}}
\subfigure[]{\includegraphics[width=4.5cm]{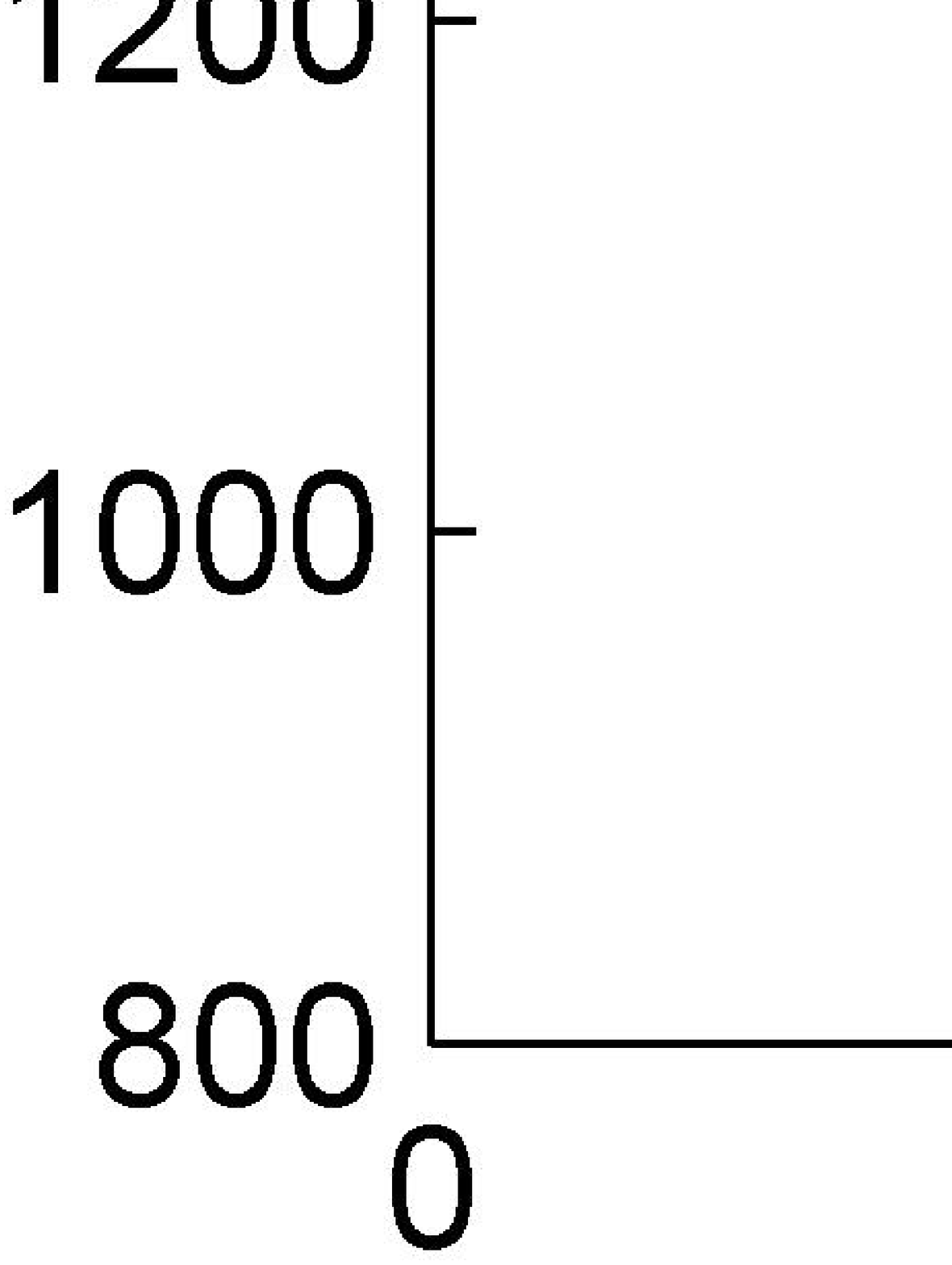}}
\captionsetup{width=14.1cm}
\caption[]{\textsl{V(1,0)-preconditioned Krylov method iteration count
    for $\sigma = -64000$ as a function of the complex shift $\beta$ for
    the 2D model problem with $N=256$. The applied method is GMRES
    with $\mu=3$ (a), $\mu=5$ (b) and $\mu=10$ (c). Minima can be
    found at $\beta = 0.26$ (a), $\beta = 0.34$ (b) and $\beta = 0.34$
    (c).}}
\label{fig:2DGMRES256}
\end{figure}

\subsection{Minimality of the complex shift w.r.t.~multigrid convergence}

To validate the results from the previous section, we use a
multigrid algorithm on a Complex Shifted Laplacian preconditioner for some well-known Krylov
methods (see \cite{erlangga2006novel}), which we apply to the
discretized model problem (\ref{eq:mod1}) with a standard all-one
right-hand side $f = 1$. This leads to an outer Krylov iteration,
where in each preconditioning phase a number of $\mu$ two-grid or V-cycle inner
steps are applied to solve the preconditioning system. Note that all experiments 
presented in this chapter are performed using left preconditioning. Two standard
Krylov methods for solving this model problem considered here are
GMRES \cite{saad1986gmres} and BiCGStab \cite{van1992bicg}. The number
of multigrid steps per preconditioning phase is denoted $\mu$, where e.g.~GMRES with $\mu=5$ 
indicates five inner multigrid iterations are used
within every outer GMRES step. As initial guess for the solution, we use a
standard all-zero vector $\bold{v}_0$.

A first result is displayed in Figure \ref{fig:1DBiCGTG}, where the
number of Krylov iterations is plotted for different choices of
$\mu$, the number of two-grid iterations per Krylov step. A relative reduction tolerance
of $10^{-6}$ is used on the initial residual and the maximum number of Krylov iterations is
capped at $64$, being the number of unknowns in the system. 
For most values of $\beta < \beta_{\min}$ the number of Krylov 
iterations required to reach the solution is excessively large. 
Note that for some $\beta < \beta_{\min}$ the destructive effect of the divergent modes 
on the global convergence appears only when sufficiently many two-grid
iterations are applied. This is due to the fact that divergence is a
limit concept, implying the error may not increase during the initial
$m_0$ multigrid iterations. However, for all shifts $\beta < \beta_{\min}$ 
the number of Krylov iterations increases dramatically as $\mu$ grows 
larger, in which case we approach the theoretical curve. 
These experimental results confirm that choosing the minimal complex 
shift at least as large as $\beta_{\min}$ always ensures a safe choice for $\beta$, 
independently of the number of multigrid cycles performed.
The experiment can easily be extended to a full V(1,0)-cycle in 1D and
2D as displayed by Figure \ref{fig:1D2DBiCGV}, supporting the
theoretical results shown on Figure \ref{fig:1Dbeta34grid} and
\ref{fig:2Dbeta34grid} for the higher-order $k$-grid schemes.

We can state that $\beta_{\min}$ is indeed \emph{minimal} with respect to multigrid
convergence, meaning it can be considered a lower limit to \emph{guarantee} convergence. 
Consequently, for a large (infinite) number of
multigrid applications, $\beta_{\min}$ is the smallest possible
complex shift for any multigrid preconditioned Krylov method to
converge.

\subsection{Near-optimality of the complex shift w.r.t.~Krylov convergence}

Another important observation can be made, relating the minimal complex shift parameter $\beta_{\min}$ to Krylov convergence behaviour. For practical purposes, we consider the slightly smaller $N=32$ model problem, with theoretical 4-grid $\beta_{\min}$-curve as shown by Figure \ref{fig:2DbetaN32}(a). 

Temporarily focussing on a specific wavenumber $\sigma$, Figure \ref{fig:2DGMRESTGIters} shows the number of Krylov iterations as a function of the complex shift $\beta$ for four different numbers of multigrid applications. Note that for $\beta \ll \beta_{\min}$ the number of Krylov iterations is typically large as discussed in the previous section. Additionally however, one clearly observes the existence of a minimum number of Krylov iterations corresponding to a certain complex shift. It can be derived from the figure that for $\sigma = -1000$ the minimum is reached at $\beta = 0.30$ (a), $\beta = 0.36$ (b), $\beta = 0.40$ (c) and $\beta = 0.42$ (d) for $\mu = 1$ (a) , $\mu = 3$ (b), $\mu = 5$ (c) and $\mu = 10$ (d) two-grid iterations respectively. These shifts are \emph{optimal} for the Krylov convergence, in the sense that they reduce the global number of iterations to a minimum. Note that the theoretical minimal complex shift for $\sigma=-1000$ equals $\beta_{\min} \approx 0.42$. Hence, it appears that for a large number of multigrid iterations, the `iteration-minimum-$\beta$' approximates the theoretical $\beta_{\min}$.

A similar tendency applies to other values of the wavenumber, as shown by Figure \ref{fig:2DGMRESTGItersLowest}, where the complex shift corresponding to the minimum amount of Krylov iterations is plotted as a function of the wavenumber $\sigma$. Comparing this `iteration-minimum-$\beta$' to the theoretical $\beta_{\min}$-curve, one observes that the latter is being approximated by the iteration-minimum curves for increasing numbers of multigrid applications. Consequently, we call $\beta_{\min}$ \emph{near-optimal} with respect to the Krylov convergence, implying that for sufficiently large numbers of multigrid iterations, the global number of Krylov iterations corresponding to $\beta_{\min}$ will be minimal. However, when approximately solving the preconditioning problem using only a small number of multigrid cycles, as is common practice, the iteration-minimum-$\beta$ may be slightly smaller than $\beta_{\min}$, hence the term `near-optimality'. Nonetheless, choosing $\beta = \beta_{\min}$ as the complex shift provides a $100\%$ multigrid stable and low-Krylov iteration solution method.

As an extra validation of these results, Figure \ref{fig:2DGMRES256} shows the number of Krylov iterations as a function of $\beta$ for a more realistic $256\times256$-grid 2D model problem preconditioned by different amounts of V(1,0)-cycles. Again, the near-optimality of $\beta_{\min}$ is clearly visible, as the shifts minimizing the number of Krylov iterations $\beta = 0.26$ (a), $\beta = 0.34$ (b) and $\beta = 0.34$ (c) approximate the theoretical minimal complex shift parameter $\beta \approx 0.34$.

\subsection{On the notion of near-optimality}

As described in the previous section, the Local Fourier Analysis from Section 2 does not always yield the most optimal shift parameter $\beta$ for a given multigrid preconditioned Krylov solver. Indeed, for some problems choosing $\beta < \beta_{\min}$ leads to faster Krylov convergence, despite the slightly diverging multigrid preconditioner. The analysis from section 2 can however be extended to the full preconditioned-Krylov iteration matrix as presented in \cite{wienands2000fourier}. For a preconditioning $k$-grid scheme, the eigenmatrix containing the Fourier symbols corresponding to the full iteration matrix $K_1^k$  can be written as 
\begin{equation} \label{eq:wien}
\tilde{K}_1^k = \tilde{A}_1(\sigma) (I_1 -  \tilde{M}_1^k) \tilde{A}_1^{-1}(\tilde\sigma),
\end{equation}
where $\tilde{A}_1(\sigma)$ is the (diagonal) eigenmatrix corresponding to the fine-grid discretization matrix $A_1(\sigma)$ of the original Helmholtz problem, $\tilde{A}_1(\tilde\sigma)$ is the Fourier representation of the perturbed preconditioning matrix $A_1(\tilde\sigma)$ and $\tilde{M}_1^k$ is the $k$-grid eigenmatrix as presented in Section 2. For completeness we note that the above equation corresponds to a MG-Krylov scheme where right-preconditioning is used.

It is shown in \cite{wienands2000fourier} that, using (\ref{eq:wien}), a sharp theoretical upper bound on the total convergence factor $\rho_i$ can be determined when applying MG-GMRES(m) to the Poisson problem. In the setting of the Helmholtz problem, the determination of a sharp upper bound for $\rho_i$ will appear to be rather non-trivial (see further). However, a sufficient condition for a given multigrid preconditioned Krylov method to converge is the so-called \emph{half-plane condition} (HPC), which states that the field of values of the iteration matrix $\{\bar{x}^T K_1^k x / \bar{x}^T x : x \in \mathbb{C}^{(N-1)^d}\}$ should be contained in an open half-plane $\{z : Re(e^{-i\varphi}z) > 0\}$ for some $\varphi \in [0,2\pi]$. This condition can be directly transferred into the LFA setting by demanding that the eigenvalues of $\tilde{K}_1^k$ lie within an open half-plane.

\begin{figure}[t] \centering
\includegraphics[width=7cm]{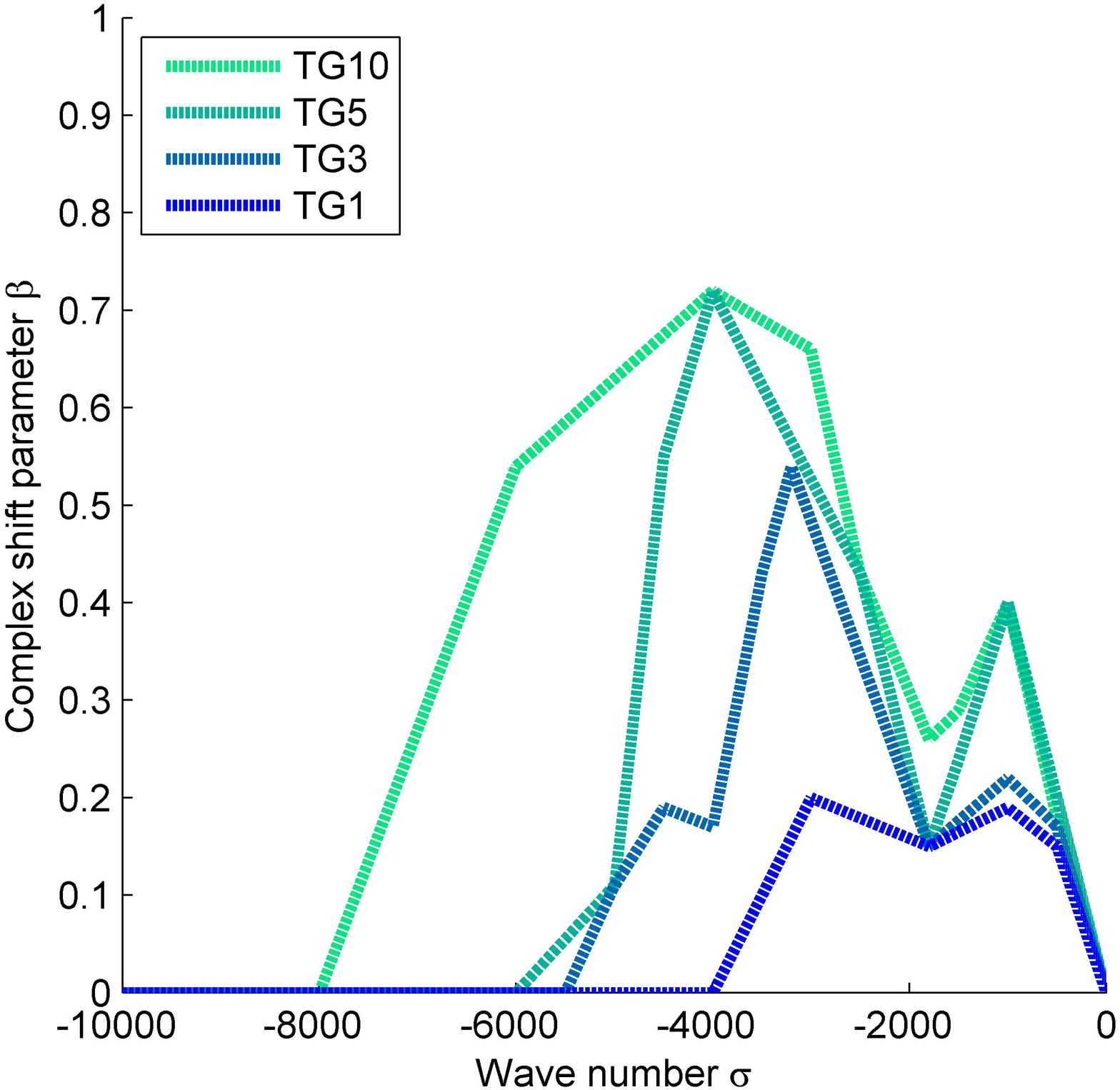}
\captionsetup{width=14.1cm}
\caption[2D two-grid preconditioned krylov half-plane condition (HPC)]{\textsl{Two-grid-preconditioned Krylov method half-plane condition minimum $\beta$ as a function of the wavenumber $\sigma$ for the 2D model problem with $N=32$. Results given for $\mu=1$, $\mu=3$, $\mu=5$ and $\mu=10$.}}
\label{fig:HPC}
\end{figure}

Figure \ref{fig:HPC} now shows a new minimal $\beta$-curve based upon the HPC convergence criterion, effectively displaying the smallest possible value of $\beta$ for which the half-plane condition is satisfied and yielding an absolute lower limit on the shift for two-grid preconditioned Krylov convergence in function of the wavenumber $\sigma$. Note that the lower limit for $\beta$ displayed on Figure \ref{fig:HPC} is distincly smaller than the corresponding iteration-minimum-$\beta$ derived from Figure \ref{fig:2DGMRESTGItersLowest}. Consequently, the shift parameter minimizing the general number of Krylov iterations can in practice be found in between the lower bound given by the HPC condition (implying MG-Krylov convergence) and the upper bound given by the LFA analysis in Section 2 (ensuring multigrid convergence), tending more towards the latter one as $\mu$ grows.

Aiming to rigourously define the iteration-minimum-$\beta$ for a general multigrid preconditioned Krylov method, the subsequent analysis from \cite{wienands2000fourier} can be transferred to a Helmholtz setting to obtain an upper bound for the convergence factor $\rho_i$. Indeed, by containing the spectrum of $K_1^k$ within an ellipse $E(c,d,a)$ (excluding the origin) with center $c$, focal distance $d$ and major-semi axis $a$, it is shown in \cite{wienands2000fourier} that $\rho_{i\gg1}$ can be heuristically estimated by
\begin{equation} \label{eq:wien2}
\rho_i \approx \frac{a+\sqrt{a^2-d^2}}{c+\sqrt{c^2-d^2}}.
\end{equation}
This estimate is generally sharp for the Poisson case as the interior of $E(c,d,a)$ is well covered by the spectrum of $K_1^k$. For most Helmholtz problems however, the ellipse constructed is rather pathological as the spectrum of $K_1^k$ only covers a fraction of $E(c,d,a)$ due to its highly irregular shape. Consequently, estimate (\ref{eq:wien2}) given above is generally \emph{not} sharp for Helmholtz problems, as illustrated by Table \ref{tab:convfac}. One observes that the iteration-minimum-$\beta$ value is clearly reflected in the experimental convergence factor estimates but is not captured by the theoretical estimates. Concluding, the ellipse-fitting methodology described in \cite{wienands2000fourier} unfortunately does not yield appropriate results for general Helmholtz problems. Future work may include the derivation of a sharp theoretical convergence factor estimate for general Helmholtz problems, leading to a rigourous prediction of the iteration-minimum-$\beta$.

\begin{table}[t!]
\centering
\begin{tabular}{c c c c c c}
\hline
  $\beta$       & 0.1 & 0.2 & 0.3 & 0.4 & 0.5 \\
\hline
  $\rho_{TH}$  		  & 1.00 & 1.00 & 0.99 & 0.99 & 0.99 \\
  $\rho_{EX}$ 			& 0.78 & 0.75 & 0.74 & 0.75 & 0.76 \\
\hline \\
\end{tabular}
\caption{\textsl{Experimental and predicted two-grid preconditioned Krylov convergence factors for the 2D model problem with $N=32$, $\mu = 1$ and wavenumber $\sigma = -1000$. Experimental factors calculated based upon Figure \ref{fig:2DGMRESTGIters}; theoretical factors based upon estimate (\ref{eq:wien2}).}}
\label{tab:convfac}
\end{table}
 
\subsection{General remarks on the number of Krylov iterations of a stable MG-Krylov solver}

This conclusive section provides the reader with some general intuition regarding the number of Krylov iterations for a sufficiently shifted (thus stable) multigrid preconditioned Krylov problem, and the dependency of this number of iterations on the wavenumber $\sigma$. As demonstrated by experimental Figures \ref{fig:1DBiCGTG}, \ref{fig:1D2DBiCGV}, \ref{fig:2DbetaN32}(b) and \ref{fig:2DGMRESTGIters}, the general number of Krylov iterations required to accurately solve the 2D model problem with a fixed and sufficiently large complex shift $\beta$ highly depends on the wavenumber $\sigma$. For sufficiently large and fixed complex shifts, one distinctly perceives the number of iterations to gradually rise to a maximum around $\sigma = -4/h^2$. Upon reaching this maximum the number of iterations decreases slowly, exhibiting a steep descend for values of $\sigma$ around $-8/h^2$. The same observation has been reported in \cite{reps2012analyzing}, where an eigenvalue analysis was performed to rigourously anticipate the Krylov convergence behaviour. Note that we have used standard Dirichlet boundary conditions, as opposed to the absorbing boundary conditions used in \cite{reps2012analyzing}; however, the conclusions are identical. The large difference in the number of Krylov iterations required is due to the indefinite nature of the problem, reaching a maximum around $\sigma = -4/h^2$ (where the problem is heavily indefinite and thus hard to solve) and causing the number of iterations to suddenly drop for values of $\sigma$ around $-8/h^2$, i.e.~where the 2D problem turns negative definite (and thus again easy to solve), as Figure \ref{fig:2DbetaN32}(b) clearly illustrates (see colors). A rigourous explanation for this behaviour is beyond the scope of this text. For additional information and a thorough analysis on the subject, we cordially refer the reader to \cite{reps2012analyzing}, \cite{simoncini2007recent}.

\section{Conclusions and discussion}

In this paper we have analyzed the convergence of a complex Shifted 
Laplacian preconditioned Krylov solver for the Helmholtz problem. 
A multigrid method is used to solve the Shifted Laplacian 
preconditioning system, which is a Helmholtz problem where the 
wavenumber $k^2$ is scaled by $(1 + \beta\iota)$. This results in an 
operator $-\Delta + \tilde{\sigma}$, where 
$\tilde{\sigma}=-k^2(1+\beta \iota)$ is complex valued.
  
The asymptotic multigrid convergence rate of a two-, three- and four-grid scheme were
analyzed theoretically with the aid of LFA. It is found that the convergence rate is 
mainly determined by the amplification factor maximum which appears at a single resonance frequency. 
The resulting convergence rate depends on the grid distance $h$, the complex shift parameter 
$\beta$ and the wavenumber $k$. 
By increasing the complex shift $\beta$, the maximum at the resonance decreases and the convergence rate
improves. In general, the larger $\beta$, the better the multigrid
convergence rate.

However, the larger we choose the value of $\beta$, the further the complex shifted problem
deviates from the original Helmholtz problem and the worse it will
perform as a preconditioner. Indeed, a balance needs to be
found between fast Krylov and multigrid convergence.

From the expression of the convergence rate it is possible to define
$\beta_{\min}$ as the smallest value of $\beta$ for which the
multilevel solution method of the shifted problem is stable. If
$\beta$ is taken smaller than $\beta_{\min}$, unstable modes are bound to destroy the multigrid solver.

When solving the preconditioner problem exactly using multigrid (up to discretization error order), the complex shift should always be
taken larger than $\beta_{\min}$. This ensures multigrid will converge when applied
to the shifted problem. Experiments show that in this case the
choice of $\beta=\beta_{\min}$ is optimal, leading to a minimal number of Krylov
iterations. However, when solving the preconditioner problem approximately 
using only a limited number of V-cycles, experimental results show that a minimum 
number of Krylov iterations is reached when choosing the shift parameter $\beta$ 
slightly smaller than $\beta_{\min}$. We conclude that $\beta_{\min}$ is a safe and
near-optimal choice for the complex shift parameter, ensuring multigrid stability and solving 
the problem using a (nearly) minimal number of Krylov iterations.

Additionally, we have shown that $\beta_{\min}$ depends in an irregular way on the
wavenumber $k$ and the grid distance $h$, and that it is furthermore dependent on the number of pre- and
post-smoothing steps. Choosing $\beta$ around 0.5, as is
common practice in the literature, ensures that the standard V(1,1)-cycle convergences
for all $k h \leq 0.625$, since this $\beta$ is distinctly larger than all corresponding $\beta_{\min}$. 
This indeed legitimizes the choice of $\beta = 0.5$ for practical purposes.

For problems with regional space-dependent wavenumbers $k(x)$, each region 
can be associated with a certain minimal complex shift $\beta_{\min}$. 
In that case, an indisputably safe choice for the complex shift parameter $\beta$
would be the largest possible regional $\beta_{\min}$ appearing in the problem. 

Note that we have used homogeneous Dirichlet boundary conditions throughout the
paper, while many applications use absorbing boundary conditions. This
shifts the eigenvalues of the original Helmholtz problem away from the
real axis and generally leads to better Krylov convergence. Consequently, the 
analysis in this paper is a discussion of the worst-case convergence scenario, 
and in practice much better Krylov convergence will be found than shown by the 
presented experiments. As Local Fourier Analysis makes no restrictions on the 
boundaries however, the theoretical results presented within this text are generally valid.

\newpage

\section*{Acknowledgements}
The authors acknowledge financial support from FWO-Flanders through project
G.0.120.08. W.V. is also supported by FWO \textit{krediet aan navorser}
project number 1.5.145.10.  Additionally, we would like to thank Bram Reps for fruitful
discussions on the subject.

\nocite{*}
\bibliographystyle{wileyj}
\bibliography{refs}

\end{document}